%%%%%%%%%%%%%%%%%%%%%%%%%%%%%%%%%%%%%%%%%%%%%%%%%%%%%%%%%%
%%
%%     This is the AMS-LaTeX file:
%%
%%     Colli-Gilardi-Sprekels 
%%     Cahn-Hilliard dynamic boundary conditions
%%
%%%%%%%%%%%%%%%%%%%%%%%%%%%%%%%%%%%%%%%%%%%%%%%%%%%%%%%%%%
%
%%%%%%%%%%%%%%%%%%%%%%%%%%%%%%%%%%%%%%%%%%%%%%%%%%%%%%%%%%
%
%   IL TESTO DEL LAVORO HA DUE VERSIONI: QUELLA SOPRA,
%   CHE E' LA PIU` RECENTE, IN EFFETTI RIVEDUTA E 
%   CORRETTA DALLA PRIMA, PREPARATA PER LA RIVISTA --
%   QUELLA SOTTO, SEMPRE IN LATEX, E' QUELLA  
%   UTILIZZATA PER IL PREPRINT ARXIV
%
%%%%%%%%%%%%%%%%%%%%%%%%%%%%%%%%%%%%%%%%%%%%%%%%%%%%%%%%%%%%%%%%%%
%%%%%%%%%%%%%%%%%%%%%%%%%%%%%%%%%%%%%%%%%%%%%%%%%%%%%%%%%%%%%%%%%%
%%%%%%%%%%%%%%%%%%% VERSIONE RIVISTA %%%%%%%%%%%%%%%%%%%%%%%%%%%%%
%%%%%%%%%%%%%%%%%%%%%%%%%%%%%%%%%%%%%%%%%%%%%%%%%%%%%%%%%%%%%%%%%%
%%%%%%%%%%%%%%%%%%%%%%%%%%%%%%%%%%%%%%%%%%%%%%%%%%%%%%%%%%%%%%%%%%

\def\LaTeX{%
  \let\Begin\begin
  \let\End\end
  \let\salta\relax
  \let\finqui\relax
  \let\futuro\relax}

\def\UK{\def\our{our}\let\sz s}
\def\USA{\def\our{or}\let\sz z}

\UK
%\USA

%%%%%%%%%%%%%%%%%%%%%%%%%%%%%%%%%

% scegliere fra \TeX e \LaTeX  e fra  \UK oppure \USA

%\TeX
\LaTeX

%\UK
\USA

%%%%%%%%%%%%%%%%%%%%%%%%%%%%%%%%%
%% page layout
%%%%%%%%%%%%%%%%%%%%%%%%%%%%%%%%%

\salta

\documentclass[twoside,12pt]{article}
\setlength{\textheight}{24cm}
\setlength{\textwidth}{16cm}
\setlength{\oddsidemargin}{2mm}
\setlength{\evensidemargin}{2mm}
\setlength{\topmargin}{-15mm}
\parskip2mm

%%%%%%%%%%%%%%%%%%%%%%%%%%%%%%%%%
%% packages
%%%%%%%%%%%%%%%%%%%%%%%%%%%%%%%%%

\usepackage{color}
\usepackage{amsmath}
\usepackage{amsthm}
\usepackage{amssymb}
\usepackage[mathcal]{euscript}

%\usepackage[notref,notcite]{showkeys}
%\usepackage{showkeys}
%
%		COLORS FOR CORRECTIONS
%
% do the same, please (i.e., don't use the standard {\color{red} text} or similar): 
% just replace dots by the color you prefer and cancel `%' in front of \def\yourname

% example of use:  \gianni{I want this to become red}

\definecolor{viola}{rgb}{0.3,0,0.7}
\definecolor{ciclamino}{rgb}{0.5,0,0.5}

\def\pier #1{{\color{red}#1}}

\def\pier #1{#1}

%%%%%%%%%%%%%%%%%%%%%%%%%%%%%%%%%
%% bibliographystyle
%%%%%%%%%%%%%%%%%%%%%%%%%%%%%%%%%

\bibliographystyle{plain}

%%%%%%%%%%%%%%%%%%%%%%%%%%%%%%%%%
%% environments
%%%%%%%%%%%%%%%%%%%%%%%%%%%%%%%%%

%

\finqui

\def\Beq{\Begin{equation}}
\def\Eeq{\End{equation}}
\def\Bsist{\Begin{eqnarray}}
\def\Esist{\End{eqnarray}}

\def\Bthm{\Begin{theorem}}
\def\Ethm{\End{theorem}}
\def\Blem{\Begin{lemma}}
\def\Elem{\End{lemma}}

\def\Bcor{\Begin{corollary}}
\def\Ecor{\End{corollary}}
\def\Brem{\Begin{remark}\rm}
\def\Erem{\End{remark}}

\def\Bdim{\Begin{proof}}
\def\Edim{\End{proof}}
\def\Bcenter{\Begin{center}}
\def\Ecenter{\End{center}}
\let\non\nonumber

%%%%%%%%%%%%%%%%%%%%%%%%%%%%%%%%%
%% macros
%%%%%%%%%%%%%%%%%%%%%%%%%%%%%%%%%

% macro salvate

% sottosezioni non numerate

\def\step #1 \par{\medskip\noindent{\bf #1.}\quad}

% abbreviazioni di parole

\def\Lip{Lip\-schitz}
\def\Holder{H\"older}
\def\aand{\quad\hbox{and}\quad}

\def\lhs{left-hand side}
\def\rhs{right-hand side}
\def\sfw{straightforward}

% versioni inglesi (UK) o americane (USA)

\def\generaliz{generali\sz}
\def\lineariz{lineari\sz}

\def\regulariz{regulari\sz}

% bold, cal e mathop

\def\multibold #1{\def\arg{#1}%
  \ifx\arg\pto \let\next\relax
  \else
  \def\next{\expandafter
    \def\csname #1#1#1\endcsname{{\bf #1}}%
    \multibold}%
  \fi \next}

\def\pto{.}

\def\multical #1{\def\arg{#1}%
  \ifx\arg\pto \let\next\relax
  \else
  \def\next{\expandafter
    \def\csname cal#1\endcsname{{\cal #1}}%
    \multical}%
  \fi \next}

% operatori

\def\multimathop #1 {\def\arg{#1}%
  \ifx\arg\pto \let\next\relax
  \else
  \def\next{\expandafter
    \def\csname #1\endcsname{\mathop{\rm #1}\nolimits}%
    \multimathop}%
  \fi \next}

\multibold
qwertyuiopasdfghjklzxcvbnmQWERTYUIOPASDFGHJKLZXCVBNM.

\multical
QWERTYUIOPASDFGHJKLZXCVBNM.

\multimathop
dist div dom meas sign supp .

% accorpamenti di formule citate:
% uso  \accorpa {prima}{seconda}
%      \Accorpa\cs prima seconda (con il comodo blank anche dopo)
% NB: \Accorpa definisce \cs come l'accorpamento delle due citazioni
% e scrive sul file.log

\def\accorpa #1#2{\eqref{#1}--\eqref{#2}}
\def\Accorpa #1#2 #3 {\gdef #1{\eqref{#2}--\eqref{#3}}%
  \wlog{}\wlog{\string #1 -> #2 - #3}\wlog{}}

% macro comode

\def\Neto{\mathrel{{\scriptstyle\nearrow}}}
\def\Seto{\mathrel{{\scriptstyle\searrow}}}

\def\graffe #1{\mathopen\{#1\mathclose\}}

\def\<#1>{\mathopen\langle #1\mathclose\rangle}
\def\norma #1{\mathopen \| #1\mathclose \|}

\def\iot {\int_0^t}
\def\ioT {\int_0^T}
\def\intQt{\int_{Q_t}}
\def\intQ{\int_Q}
\def\iO{\int_\Omega}
\def\iG{\int_\Gamma}
\def\intS{\int_\Sigma}
\def\intSt{\int_{\Sigma_t}}

\def\dt{\partial_t}
\def\dn{\partial_n}

\def\cpto{\,\cdot\,}

\def\checkmmode #1{\relax\ifmmode\hbox{#1}\else{#1}\fi}
\def\aeO{\checkmmode{a.e.\ in~$\Omega$}}
\def\aeQ{\checkmmode{a.e.\ in~$Q$}}
\def\aeG{\checkmmode{a.e.\ on~$\Gamma$}}
\def\aeS{\checkmmode{a.e.\ on~$\Sigma$}}
\def\aet{\checkmmode{a.e.\ in~$(0,T)$}}

\def\aaQ{\checkmmode{for a.a.~$(x,t)\in Q$}}

\def\aat{\checkmmode{for a.a.~$t\in(0,T)$}}

% insiemi numerici

\def\erre{{\mathbb{R}}}

% spazi di funzioni a valori vettoriali su [0,T], [0,t], [0,s], [0,+\infty), [\delta,T]

% Come ricordare: in generale i simboli L H W  C da soli per gli spazi su (0,T)
% gli stessi raddoppiati per (0,+\infty)
% aggiunta di t o s al simbolo per (0,t) e (0,s)
% aggiunta di d al simbolo semplice o doppio per intervalli (\delta,T) e (\delta,+\infty)
% il simbolo C e i suoi derivati mettono le quadre anziche' le tonde

% Esempi   \L2V   \L\infty\Vp   \W{1,1}H   \C0H   \LL2V   \CC0\Vp   \Ld2V  \CCdH

\def\genspazio #1#2#3#4#5{#1^{#2}(#5,#4;#3)}
\def\spazio #1#2#3{\genspazio {#1}{#2}{#3}T0}

\def\L {\spazio L}
\def\H {\spazio H}
\def\W {\spazio W}

\def\C #1#2{C^{#1}([0,T];#2)}

% spazi di funzioni su \Omega, \Gamma, Q e \Sigma

\def\Lx #1{L^{#1}(\Omega)}
\def\Hx #1{H^{#1}(\Omega)}
\def\Wx #1{W^{#1}(\Omega)}
\def\LxG #1{L^{#1}(\Gamma)}
\def\HxG #1{H^{#1}(\Gamma)}

\def\LQ #1{L^{#1}(Q)}
\def\LS #1{L^{#1}(\Sigma)}

\def\Luno{\Lx 1}
\def\Ldue{\Lx 2}
\def\Linfty{\Lx\infty}
\def\LinftyG{\LxG\infty}
\def\Huno{\Hx 1}
\def\Hdue{\Hx 2}

\def\HunoG{\HxG 1}
\def\HdueG{\HxG 2}

\def\LunoG{\LxG 1}
\def\LdueG{\LxG 2}

% spazi di funzioni su Q e S

\def\LQ #1{L^{#1}(Q)}

\def\LS #1{L^{#1}(\Sigma)}

% lettere greche

\let\theta\vartheta
\let\eps\varepsilon
\let\lam\lambda

\let\TeXchi\chi                         % new \chi, exactly on the baseline
\newbox\chibox
\setbox0 \hbox{\mathsurround0pt $\TeXchi$}
\setbox\chibox \hbox{\raise\dp0 \box 0 }
\def\chi{\copy\chibox}

% quadratino di fine dimostrazione

\def\QED{\hfill $\square$}

% abbreviazioni specifiche del lavoro

\def\suG{{\vrule height 5pt depth 4pt\,}_\Gamma}

\def\yG{y_\Gamma}
\def\xiG{\xi_\Gamma}
\def\yz{y_0}
\def\xiz{\xi_0}
\def\xiGz{\xi_{\Gamma\!,0}}
\def\ustar{u^*}
\def\vstar{v^*}

\def\VG{V_\Gamma}
\def\HG{H_\Gamma}
\def\Vp{V^*}

\def\normaV #1{\norma{#1}_V}
\def\normaH #1{\norma{#1}_H}

\def\normaHG #1{\norma{#1}_{\HG}}
\def\normaVp #1{\norma{#1}_*}

\def\nablaG{\nabla_{\!\Gamma}}
\def\DeltaG{\Delta_\Gamma}

\let\hat\widehat
\def\Beta{\hat{\vphantom t\smash\beta\mskip2mu}\mskip-1mu}
\def\betaeps{\beta_\eps}
\def\Betaeps{\Beta_\eps}
\def\betaG{\beta_\Gamma}
\def\BetaG{\Beta_\Gamma}
\def\betaz{\beta^\circ}
\def\betaGz{\betaG^\circ}

\def\betaGeps{\beta_{\Gamma\!,\eps}}
\def\BetaGeps{\Beta_{\Gamma\!,\eps}}
\def\piG{\pi_\Gamma}
\def\Pi{\hat\pi}
\def\PiG{\Pi_\Gamma}
\def\gG{g_\Gamma}
\def\lamG{\lam_\Gamma}

\def\feps{f_\eps}
\def\fGeps{f_{\Gamma\!,\eps}}

\def\mz{m_0}
\def\Bp{B_{p,\eps}}

\def\yeps{y_\eps}
\def\yGeps{y_{\Gamma\!,\eps}}
\def\weps{w_\eps}

\def\xieps{\xi_\eps}
\def\xiGeps{\xi_{\Gamma\!,\eps}}

\def\taueps{\tau_\eps}
\def\geps{g_\eps}

% scelte che si possono cambiare

\def\cd{c_\delta}
\def\deltaz{\delta_0}
\def\Cz{C_0}

%%%%%%%%%%%%%%%%%%%%%%%%%%%%%%
\Begin{document}
%%%%%%%%%%%%%%%%%%%%%%%%%%%%%%%%%

%%%%%%%%%%%%%%%%%%%%%%%%%%%%%%%%%
%% front page
%%%%%%%%%%%%%%%%%%%%%%%%%%%%%%%%%

\title{On the Cahn\pier{--}Hilliard equation\\[0.3cm]
  with dynamic boundary conditions\\[0.3cm]
  and a dominating boundary potential
\footnote{{\bf Acknowledgment.}\quad\rm
The present note benefits from a partial support of the MIUR-PRIN Grant 2010A2TFX2 ``Calculus of variations'' and the GNAMPA (Gruppo Nazionale per l'Analisi Matematica, la Probabilit\`a e le loro Applicazioni) of INdAM (Istituto Nazionale di Alta Matematica) for PC and GG.}}
\author{}
\date{}
\maketitle
\Bcenter
\vskip-1cm
{\large\sc Pierluigi Colli$^{(1)}$}\\
{\normalsize e-mail: {\tt pierluigi.colli@unipv.it}}\\[.25cm]
{\large\sc Gianni Gilardi$^{(1)}$}\\
{\normalsize e-mail: {\tt gianni.gilardi@unipv.it}}\\[.25cm]
{\large\sc J\"urgen Sprekels$^{(2)}$}\\
{\normalsize e-mail: {\tt sprekels@wias-berlin.de}}\\[.45cm]
$^{(1)}$
{\small Dipartimento di Matematica ``F. Casorati'', Universit\`a di Pavia}\\
{\small via Ferrata 1, 27100 Pavia, Italy}\\[.2cm]
$^{(2)}$
{\small Weierstra\ss-Institut f\"ur Angewandte Analysis und Stochastik}\\
{\small Mohrenstra\ss e\ 39, 10117 Berlin, Germany}\\[.8cm]
\Ecenter

\Begin{abstract}
\pier{The Cahn\pier{--}Hilliard and viscous Cahn\pier{--}Hilliard equations with singular
and possibly nonsmooth potentials and dynamic boundary condition 
are considered and some well-posedness and regularity results are proved.}
\vskip3mm

\End{abstract}

\salta

\pagestyle{myheadings}
\newcommand\testopari{\sc Colli \ --- \ Gilardi \ --- \ Sprekels}
\newcommand\testodispari{\sc Cahn\pier{--}Hilliard equation
with dynamic boundary conditions}
\markboth{\testodispari}{\testopari}

\finqui

%%%%%%%%%%%%%%%%%%%%%%%%%%%%%%%%%
%% very beginning
%%%%%%%%%%%%%%%%%%%%%%%%%%%%%%%%%

\section{Introduction}
\label{Intro}
\setcounter{equation}{0}

The classical Cahn\pier{--}Hilliard equation and the so-called viscous Cahn\pier{--}Hilliard equation (see \cite{CahH, EllSt, EZ}) in their simplest forms read 
\Beq
  \dt y - \Delta w = 0
  \aand
  w = \tau \dt y - \Delta y + \beta(y) + \pi(y) - g
  \quad \hbox{in $\Omega\times (0,T) $},
  \label{cahnh}
\Eeq
according to the case $\tau =0$ or $\tau >0,$ respectively. Some significant extensions and a comparative discussion on the modeling approach for phase separation and atom mobility between cells can be found in \cite{FG, Gu, Podio, CGPS3}).

In~\eqref{cahnh}, $y$ denotes the order parameter and $w $ represents 
the chemical potential.  Moreover, 
$\beta$ and $\pi$ are the derivatives of the convex part $\Beta$
and of the concave perturbation $\Pi$ of a double well potential~$\calW$,
and $g$ is a source term.
Important examples of $\calW$ are the everywhere defined regular potential $\calW_{reg}$ and the logarithmic double-well potential $\calW_{log}$ given~by
\begin{align}
  & \calW_{reg}(r) = \frac14(r^2-1)^2 \,,
  \quad r \in \erre 
  \label{regpot}
  \\
  & \calW_{log}(r) = ((1+r)\ln (1+r)+(1-r)\ln (1-r)) - c r^2 \,,
  \quad r \in (-1,1)
  \label{logpot}
\end{align} 
where $c>0$ in the latter is large enough in order to kill convexity. Another important example refers to the so-called double-obstacle problem and corresponds to the nonsmooth potential ${\cal W}_{2obst}$ specified by
\Beq
{\cal W}_{2obst}(r)= \left\{ 
\begin{array}{cl}
c(1-r^2) & \hbox{ if } \,|r|\leq 1  
\\[0.1cm]
+\,\infty  & \hbox{ if } \,|r| >1  
\end{array} \right. .
\label{pier4}
\Eeq
In this case $\beta$ is no longer a derivative, but it represents the
subdifferential $\partial I_{[-1, 1]}$ of the indicator function of the interval $[-1, 1]$, that is, 
\begin{equation}
\label{pier5}
s \in \partial I_{[-1, 1]} (r) 
\quad \hbox{ if and only if } \quad 
s \ \left\{
\begin{array}{ll}
\displaystyle
\leq \, 0 \   &\hbox{if } \ r=-1   
\\[0.1cm]
= \, 0 \   &\hbox{if } \ -1< r < 1  
\\[0.1cm]
\geq \, 0 \  &\hbox{if } \  r = 1  
\\[0.1cm]
\end{array}
\right. .
\end{equation}
We are interested to the coupling of \eqref{cahnh} with the usual no-flux condition for the chemical potential
\Beq
  (\dn w)\suG =0 
  \quad \hbox{on $\Gamma\times (0,T)$}
  \label{nfc}
\Eeq
and with the dynamic boundary condition
\Beq
  (\dn y)\suG + \dt\yG - \DeltaG\yG + \betaG(\yG) + \piG(\yG) = \gG
  \quad \hbox{on $\Gamma\times (0,T)$}
  \label{dbc}
\Eeq
where $\yG $ denotes the trace $ y\suG$ on the boundary $\Gamma$ of $\Omega$, $ - \DeltaG$ stands for the Laplace\pier{--}Beltrami operator on $\Gamma$, $\betaG$ and $\piG$ are nonlinearities playing the same role as $\beta$ and $\pi$ but now acting on the boundary value of the order parameter, and finally $\gG$ is a boundary source term with no relation with $g$ acting on the bulk.  

The physical meaning and free energy derivation of the boundary value problem given by \eqref{cahnh} and \eqref{nfc}--\eqref{dbc} have been discussed specifically in \cite{GiMiSchi}.  The Cahn\pier{--}Hilliard equation \eqref{cahnh}, endowed with the dynamic boundary condition \eqref{dbc}, has drawn much attention in recent years: we quote, among other contributions, \cite{CFP, Is, MirZelik, PRZ, RZ, WZ}. In particular, the existence and
uniqueness of solutions as well as the behavior of the solutions as time goes to 
infinity have been studied for regular potentials $\calW$ and $\calW_\Gamma= \BetaG + \PiG$. Moreover, a wide class of potentials, 
including especially singular potentials like \eqref{logpot} and \eqref{pier4}, 
has been considered in \cite{GiMiSchi, GiMiSchi2}: in these two papers the authors 
were able to overcome the difficulties due to singularities and to show 
well-posedness results along with the long-time behavior of solutions. The approach of \cite{GiMiSchi, GiMiSchi2} is based on a set of assumptions for $\beta, \, \pi$ 
and $\betaG, \, \piG $ that gives the role of the dominating potential to $\calW$ instead of to ${\calW}_\Gamma$ and entails some technical difficulties. 

In this note, we follow a strategy developed in \cite{CaCo} to investigate the Allen\pier{--}Cahn equation with 
dynamic boundary condition, which consists in letting ${\calW}_\Gamma$ be the 
leading potential between the two. This approach simplifies the analysis and allows for
a unified treatment of the initial value problem for \eqref{cahnh}, \eqref{nfc},  \eqref{dbc} and for a linearized version thereof. This was a main motivation for
this paper, namely to complement and improve the results of \cite{GiMiSchi}. Another input for the realization of this article was the related project of investigating the optimal control problem for the Cahn\pier{--}Hilliard equation with dynamic boundary condition. In view of the already realized contributions for the 
corresponding Allen\pier{--}Cahn equation (see \cite{CS} and \cite{CFS}), a work 
program for the more difficult Cahn\pier{--}Hilliard setting appeared to be 
natural and worth pursuing. This will be the subject of a forthcoming contribution, which will make intense use of the results established here. 

Concerning the optimal control problems, let us mention that in \cite{CS} both the cases of distributed and boundary controls have been addressed for logarithmic-type 
potentials as in \eqref{logpot}: after showing the existence of optimal controls and checking that the control-to-state mapping is twice continuously Fr\'echet 
differentiable, first-order necessary optimality conditions were established  in terms of a variational inequality and the adjoint state equation, and 
second-order sufficient optimality conditions were proved. The related paper \cite{CFS} deals with (non-differentiable) double obstacle potentials (see \eqref{pier4}) and contains the 
proofs of the existence of optimal controls and the derivation of first-order necessary conditions of optimality. Using the results from \cite{CS} 
for the case of (differentiable) logarithmic potentials, a so-called ``deep quench limit'' is performed to derive first-order necessary optimality conditions.

With the above motivation in mind, we study here the initial and boundary value problem 
\begin{align}
  & \dt y - \Delta w = 0
  \quad \hbox{in $\, Q:=\Omega\times (0,T)$}
  \label{Iprima}
  \\[0.2cm]
  & w = \tau \, \dt y - \Delta y + \beta(y) + \pi(y) - g
  \quad \hbox{in $\,Q$}
  \label{Iseconda}
  \\[0.2cm]
  & \dn w = 0
  \quad \hbox{on $\, \Sigma:= \Gamma\times (0,T)$}
  \label{Ibc}
  \\[0.2cm]
  & \yG = y\suG
  \aand
  \dt\yG + (\dn y)\suG - \DeltaG\yG + \betaG(\yG) + \piG(\yG) = \gG
  \quad \hbox{on $\, \Sigma$}
  \label{Iterza}
  \\[0.2cm]
  & y(0) = \yz
  \quad \hbox{in $\, \Omega$}
  \label{Icauchy}
\end{align}
\Accorpa\Ipbl Iprima Icauchy
as well as a linearization thereof, in which $\beta(y)+\pi(y)$ and $\betaG(\yG)+\piG(\yG)$ are replaced by $\lam\,y$ and $\lamG\yG$, for 
some given and a.e.\ bounded functions $\lam$ and $\lamG$ on $Q$ and~$\Sigma$, respectively.
We investigate both the viscous case $\tau >0$ and the pure case $\tau =0$, making the necessary distinctions and specifications. We show existence, uniqueness and regularity results, which are already introduced and made precise in the next section.
Section~\ref{UniqContDip} develops the details of the continuous dependence estimate that is also leading to uniqueness. The final Section~\ref{ExistenceRegularity} is concerned with the proofs of existence and of the various regularity results presented in our contribution: in particular, we prove the global boundedness of both $y$ and $\yG$ in our general framework for potentials $\calW, \, \calW_\Gamma$ and graphs~$\beta, \, \betaG$. 

%%%%%%%%%%%%%%%%%%%%%%%%%%%%%%%%%%%%%%%%%%%%%%%%%%%%%%%%%%%%%%%%%%%%%%%%

\section{Main results}
\label{MainResults}
\setcounter{equation}{0}

In this section, we describe the problem \pier{under study} and state our results.
As in the Introduction,
$\Omega$~is the body where the evolution takes place.
Being clear that just minor changes are needed to treat the lower-dimensional cases,
we assume $\Omega\subset\erre^3$
to~be open, bounded, connected, and smooth
and we write $|\Omega|$ for its Lebesgue measure.
Moreover, we still denote  
the boundary of~$\Omega$, the outward normal derivative, the surface gradient 
and the Laplace\pier{--}Beltrami operator 
by $\Gamma$, $\dn$, $\nablaG$ and $\DeltaG$, respectively.
Given a finite final time~$T$,
we set for convenience
\Bsist
  && Q_t := \Omega \times (0,t)
  \aand
  \Sigma_t := \Gamma \times (0,t)
  \quad \hbox{for every $t\in(0,T]$}
  \label{defQtSt}
  \\
  && Q := Q_T \,,
  \aand
  \Sigma := \Sigma_T \,.
  \label{defQS}
\Esist
\Accorpa\defQeS defQtSt defQS
Now, we make the assumptions on the structure of our system precise.  
However, besides the Cahn\pier{--}Hilliard equations with or without viscosity,
we are interested in solving the corresponding \lineariz ed problem around some solution as well.
As the latter corresponds to replace $\beta(y)+\pi(y)$ and $\betaG(\yG)+\piG(\yG)$ 
by $\lam\,y$ and $\lamG\yG$ in \eqref{Iseconda} and~\eqref{Iterza}, respectively,
where $\lam$ and $\lamG$ are some functions on $Q$ and~$\Sigma$,
we consider a problem that is slightly more general,
in~order to unify the treatment.
So, we assume that we are given structural functions $\Beta$,
$\BetaG$, $\pi$, $\piG$, two functions $\lam$, $\lamG$ and a constant $\tau$ satisfying the
conditions listed below.
\begin{align}
  & \hbox{$\Beta,\BetaG:\erre\to[0,+\infty]$
    are convex, proper, and l.s.c.\ and $\Beta(0)=\BetaG(0)=0$}\quad
  \label{hpBeta}
  \\
  & \hbox{$\pi,\piG:\erre\to\erre$ are \Lip\ continuous with $\pi(0)=\piG(0)=0$}
  \label{hppilip}
  \\
  & \lam \in \LQ\infty
  \aand
  \lamG \in \LS\infty
  \label{hplam}
  \\
  & \tau \geq 0 
  \label{hpvisc}
  \\
  & \lam \in \L\infty{\Wx{1,3}} 
  \quad \hbox{if $\tau=0$} \,.
  \label{hpdxlam}
\end{align}
\Accorpa\HPstruttura hpBeta hpdxlam
\Accorpa\HPstrutturataupos hpBeta hplam
We define the graphs $\beta$ and $\betaG$ in $\erre\times\erre$ by
\Beq
  \beta := \partial\Beta
  \aand
  \betaG := \partial\BetaG
  \label{defbeta}
\Eeq
and note that \pier{both} $\beta$ and $\betaG$ are maximal monotone
with some effective domains $D(\beta)$ and~$D(\betaG)$.
Due to~\eqref{hpBeta}, we have $\beta(0)\ni0$ and $\betaG(0)\ni0$ .
In the sequel, for any maximal monotone graph $\gamma:\erre\to2^{\erre}$,
we use the notation
(see, e.g., \cite[p.~28]{Brezis})
\Bsist
  && \hbox{$\gamma^\circ(r)$ is the element
    of $\gamma(r)$ having minimum modulus}
  \label{defbetaz}
  \\
  && \hbox{$\gamma_\eps^Y$ is the Yosida \regulariz ation
    of $\gamma$ at level $\eps$, for $\eps>0$}.
  \label{yosida}
\Esist
Moreover, we still write the symbol $\gamma$ 
(and,~e.g., $\gamma_\eps^Y$ as a particular case) 
for~the maximal monotone operator induced by $\gamma$ on the space~$\LQ2$.
For the graphs $\beta$ and $\betaG$ we assume the following compatibility condition
\begin{align}
  D(\betaG) \subseteq D(\beta)
  \aand
  |\betaz(r)| \leq \eta |\betaGz(r)| + C \non \qquad \qquad \\
  \hbox{for some $\eta,\, C>0$ and every $r\in D(\betaG)$}
  \label{hpcompatib}
\end{align}
and note that, roughly speaking, 
it is opposite to the one postulated in~\cite{GiMiSchi}.
On the contrary, condition~\eqref{hpcompatib}
is the same as the one introduced in the paper~\cite{CaCo},
which however deals with the Allen\pier{--}Cahn equation.
For those reasons, the results we obtain are completely new.
Next, in order to simplify the notation, we~set
\Bsist
  && V := \Huno, \quad
  H := \Ldue, \quad
  \HG := \LdueG , \quad
  \VG := \HunoG 
  \label{defspazi}
  \\
  && \calV := \graffe{v\in V:\ v\suG\in\VG}
  \label{defcalV}
\Esist
and endow the former spaces with their usual norms
and the latter with the graph norm.
For the norms in the generic Banach space~$X$
and in any power of it,
we write~$\norma\cpto_X$.
However, simpler symbols are used in particular cases.
For instance, for $1\leq p\leq\infty$, $\norma\cpto_p$ is the usual norm in $\Lx p$
and $\normaVp\cpto$~denotes a precise norm in~$\Vp$
introduced later~on
(see the forthcoming~\eqref{defnormaVp}).
Such a norm is equivalent to the standard dual norm.
Furthermore, the symbol $\<\cpto,\cpto>$ stands
for the duality pairing between $\Vp$ and~$V$.
In the sequel, it is understood that $H$ is embedded in $\Vp$
in the usual way, i.e., in order that
$\<u,v>=(u,v)$, the inner product of~$H$,
for every $u\in H$ and $v\in V$.
Finally, we define the \generaliz ed mean value of any $\vstar\in\Vp$ by setting
\Beq
  {\pier{\vstar_{\,\Omega}} := \frac 1 {|\Omega|} \, \< \vstar , 1 >
  \quad \hbox{for $\vstar\in\Vp$}.}
  \label{media}
\Eeq
Clearly, \eqref{media} gives the usual mean value when applied to elements
of~$H$.

At this point, we can describe our problem, which consists in the
variational formulation of system \Ipbl.
To prepare the assumptions we need for our main existence result
(the~problem is meaningful and a uniqueness type theorem holds under much weaker hypotheses, indeed),
we~give the data $g$, $\gG$, and $\yz$ satisfying 
\Bsist
  && g \in \L2H
  \aand
  \gG \in \L2\HG
  \label{hpg}
  \\
  && g \in \H1H
  \quad \hbox{if $\tau=0$}
  \label{hpgtauzero}
  \\
  && \yz \in \calV , \quad
  \pier{\Beta(\yz) \in \Luno 
  % , \quad
  % \Beta(\yz\suG) \in \LunoG
  \aand
  \BetaG (\yz\suG) \in \LxG 1}
  \qquad
  \label{hpyz}
  \\
  && \mz := (\yz)_\Omega 
  \quad \hbox{lies in the interior of $D(\betaG)$}.
  \label{hpinterno}
\Esist
\Accorpa\HPdati hpg hpinterno
Our problem consists in looking for a quintuplet $(y,\yG,w,\xi,\xiG)$ such~that
\begin{align}
  & y \in \H1\Vp \cap \L\infty V \cap \L2\Hdue
  \, \hbox{ and } \, \tau \, \dt y\in \L2H
  \label{regy}
  \\
  & \yG \in \H1\HG \cap \L\infty\VG \cap \L2\HdueG
  \label{regyG}
  \\
  & \yG(t) = y(t)\suG
  \quad \aat
  \label{tracciay}
  \\
  & w \in \L2V
  \label{regw}
  \\
  & \xi \in \L2H
  \aand
  \xi \in \beta(y)
  \quad \aeQ
  \label{regxi}
  \\
  & \xiG \in \L2\HG
  \aand
  \xiG \in \betaG(\yG)
  \quad \aeS
  \label{regxiG}
\end{align}
\Accorpa\Regsoluz regy regxiG
and satisfying \aat\ the variational equations
\begin{align}
  & \< \dt y(t) , v >
  + \iO \nabla w(t) \cdot \nabla v = 0
  \label{prima}
  \\
  \noalign{\smallskip}
  & \iO w(t) v
  = \iO \tau \, \dt y(t) \, v
  + \iG \dt\yG(t) \, v
  + \iO \nabla y(t) \cdot \nabla v
  + \iG \nablaG\yG(t) \cdot \nablaG v
  \qquad
  \non
  \\
  & \quad {}
  + \iO \bigl( \xi(t) + \lam(t) \, \pi(y(t)) - g(t) \bigr) \, v
  + \iG \bigl( \xiG(t) + \lamG(t) \, \piG(\yG(t)) - \gG(t) \bigr) \, v
  \label{seconda}
\end{align}
for every $v\in V$ and every $v\in\calV$, respectively,
and the Cauchy condition
\Beq
  y(0) = \yz \,.
  \label{cauchy}
\Eeq
\Accorpa\pbl prima cauchy
\Accorpa\Pbl regy cauchy
For simplicity,
we have \pier{used} the same symbol $v$ for both the test function and its trace on the boundary,
and \pier{we do so} in the sequel
if no misunderstanding can arise.
Moreover, we write products by~$\tau$ 
(like the pointwise value $\tau\dt \pier{y}(t)$ in \eqref{seconda}
which might be meaningless, in principle)
for the sake of conciseness, also in the sequel.
In such cases, it is understood that the product vanishes for $\tau=0$.
We note that an equivalent formulation of \accorpa{prima}{seconda}
is given~by
\begin{align}
  & \ioT \< \dt y(t) , v(t) > \, dt
  + \intQ \nabla w \cdot \nabla v = 0
  \label{intprima}
  \\
  \noalign{\smallskip}
  & \intQ wv
  = \intQ \tau \dt y \, v
  + \intS \dt\yG \, v
  + \intQ \nabla y \cdot \nabla v
  + \intS \nablaG\yG \cdot \nablaG v
  \qquad
  \non
  \\
  & \quad {}
  + \intQ \bigl( \xi + \lam \, \pi(y) - g \bigr) \, v
  + \intS \bigl( \xiG + \lamG \piG(\yG) - \gG \bigr) \, v
  \qquad
  \label{intseconda}
\end{align}
\Accorpa\intpbl intprima intseconda
for every $v\in\L2 V$ and every $v\in\L2\calV$, respectively.

\Brem
\label{Conserved}
Even though what we say is completely standard for Cahn\pier{--}Hilliard equations,
it is worth to note it.
By testing \eqref{prima} by the constant $1/|\Omega|$,
we~obtain
\Beq
  (\dt y(t))_\Omega = 0
  \quad \aat
  \aand
  y(t)_\Omega = \mz
  \quad \hbox{for every $t\in[0,T]$}
  \label{conserved}
\Eeq
with the notations~\eqref{media} and~\eqref{hpinterno}.
\Erem

As far as uniqueness and continuous dependence are concerned,
we have

\Bthm
\label{ContDip}
Assume % \HPstruttura\ and 
\pier{\accorpa{hpBeta}{defbeta}}
and let $(g_i,g_{\Gamma\!,i},y_{0,i})$, $i=1,2$, be two sets of data
satisfying \eqref{hpg} and such that 
$y_{0,1},y_{0,2}$ belong to $V$ and have the same mean value.
Then, if $(y_i,y_{\Gamma\!,i},w_i,\xi_i,\xi_{\Gamma\!,i})$
are any two corresponding solutions to problem \Pbl,
the inequality
\Bsist
  && \norma{y_1-y_2}_{\L\infty\Vp}^2
  + \tau \norma{y_1-y_2}_{\L\infty H}^2
  + \norma{y_{\Gamma\!,1}-y_{\Gamma\!,2}}_{\L\infty\HG}^2
  \non
  \\
  && \quad {}
  + \norma{\nabla(y_1-y_2)}_{\L2H}^2
  + \norma{\nablaG(y_{\Gamma\!,1}-y_{\Gamma\!,2)}}_{\L2\HG}^2
  \non
  \\
  && \leq c \,
     \Bigl\{
       \normaVp{y_{0,1}-y_{0,2}}^2
       + \tau \normaH{y_{0,1}-y_{0,2}}^2
       + \normaHG{y_{0,1}\suG-y_{0,2}\suG}^2
  \non
  \\
  && \qquad \quad {}
       + \norma{g_1-g_2}_{\L2H}^2
       + \norma{g_{\Gamma\!,1}-g_{\Gamma\!,2}}_{\L2\HG}^2
     \Bigr\}
  \label{contdip}
\Esist
holds true with a constant $c$ that depends only on
$\Omega$, $T$, the \Lip\ constants of $\pi$ and~$\piG$ 
and on the norms \pier{$\|\lam \|_{L^\infty (Q)}$ and~$\|\lamG \|_{L^\infty (\Sigma)} $.}
In particular, any two solutions to problem \Pbl\
have the same components $y$, $\yG$ and~$\xiG$.
Moreover, even the components $w$ and $\xi$ of such solutions are the same
if $\beta$ is single-valued.
\Ethm

The above theorem is quite similar to the results stated in \cite[Thm.~1 and Rem.~9]{GiMiSchi}.
In the same paper (see \cite[Rem.~4 and Rem.~8]{GiMiSchi}), it is also shown that \pier{partial uniqueness
and conditionally full uniqueness} as in the above statement are the best one can prove.
At this point, we are mainly interested in existence and regularity 
and what we prove is new with respect to~\cite{GiMiSchi}, as already observed.
Here is our first result in that direction.

\Bthm
\label{Existence}
Assume \pier{\accorpa{hpBeta}{defbeta}, 
\eqref{hpcompatib} and \HPdati}.
Then, there exists a quintuplet $(y,\yG,w,\xi,\xiG)$ satisfying \Regsoluz\
and solving problem \pbl.
\Ethm

Our next goal is regularity, and we present several results.
First, we want to prove that the unique solution to problem \pbl\
given by the above theorems also satisfies 
\begin{align}
  & y \in \W{1,\infty}\Vp \cap \H1V \cap \L\infty\Hdue
  \aand
  \tau \dt y\in \L\infty H
  \label{regybis}
  \\ 
  & \yG \in \W{1,\infty}\HG \cap \H1\VG \cap \L\infty\HdueG 
  \label{regyGbis}
\end{align}
whence also
\Beq
  y \in \LQ\infty
  \aand
  \yG \in \LS\infty .
  \label{ybdd}
\Eeq
\Accorpa\Regsoluzbis regybis ybdd
To this aim, \pier{besides \eqref{hplam} and \eqref{hpdxlam}} we suppose that 
\Beq
  \lam \in \W{1,\infty}H
  \aand
  \lamG \in \W{1,\infty}\HG \,.
  \label{hpdtlam}
\Eeq
As far as the data are concerned, we also assume
\Bsist
  && g \in \H1H
  \aand
  \gG \in \H1\HG
  \label{hpgbis}
  \\
  && \yz\in\Hdue , \quad
  \dn\yz\suG = 0
  \aand
  \yz\suG \in \HdueG
  \label{hpyzHdue}
  \\
  && \hbox{there exists $\xiz\in H$ such that $\xiz\in\beta(\yz)$ \aeO}
  \label{hpxiz}
  \\
  && \hbox{there exists $\xiGz\in\HG$ such that $\xiGz\in\betaG(\yz\suG)$ \aeG}
  \label{hpxiGz}
\Esist
\Accorpa\HPdatibis hpgbis hpxiGz
and, if $\tau=0$, we reinforce \eqref{hpxiz} by requiring that
\Beq
  \hbox{the family} \quad
  \graffe{-\Delta\yz - \betaeps(\yz) - g(0):\ \eps\in(0,\eps_0)}
  \quad \hbox{is bounded in $V$} 
  \label{hpxizV}
\Eeq
for some $\eps_0>0$.
Clearly, in order to ensure~\eqref{hpxizV}, one can assume that
$\Delta\yz+g(0)\in V$ and that $\betaeps(\yz)$ remains bounded in~$V$ for $\eps$ small enough,
and a sufficient condition for the latter is the following:
there exist $r_\pm,r_\pm'\in\erre$ such that
$r_-'<r_-\leq\yz\leq r_+<r_+'$ \aeO, $(r_-',r_+')\subset D(\beta)$
and the restriction of $\beta$ to $(r_-',r_+')$ is a single-valued \Lip\ continuous function.

Here is our first regularity result.

\Bthm
\label{Regularity}
Assume \pier{\accorpa{hpBeta}{defbeta}}, \eqref{hpcompatib} and \eqref{hpdtlam} on the structure 
and suppose that the data satisfy \HPdatibis\ and~\eqref{hpinterno}.
Then, the unique solution to problem \pbl\
given by Theorems~\ref{ContDip} and~\ref{Existence}
also satisfies~\Regsoluzbis. \pier{Moreover, we have that}
\Beq 
\pier{w \in \L\infty V, \quad 
\xi \in \L\infty H, \quad
\xiG \in \L\infty \HG .}
\label{pier0}
\Eeq
\Ethm

We present at once a consequence of our results,
which is obtained by simply taking
\Beq
  \Beta(r) = \BetaG(r) = 0 
  \aand
  \pi(r) = \piG(r) = r 
  \quad \hbox{for every $r\in\erre$} .
  \non
\Eeq

\Bcor
\label{Linearcase}
Assume $\tau>0$ and \eqref{hplam}.
Moreover, assume \eqref{hpgbis} and~\eqref{hpyzHdue}.
Then, there exists a unique triplet $(y,\yG,w)$ satisfying 
the regularity requirements \Regsoluzbis, \eqref{regw}
and the Cauchy condition~\eqref{cauchy},
and solving the variational equations \eqref{prima} for every $v\in V$ and
\Bsist
  && \iO w(t) v
  = \iO \tau \, \dt y(t) \, v
  + \iG \dt\yG(t) \, v
  + \iO \nabla y(t) \cdot \nabla v
  + \iG \nablaG\yG(t) \cdot \nablaG v
  \qquad
  \non
  \\
  && \quad {}
  + \iO \bigl( \lam(t) \, y(t) - g(t) \bigr) \, v
  + \iG \bigl( \lamG(t) \, \yG(t) - \gG(t) \bigr) \, v
  \non
\Esist
for every $v\in\calV$. 
\Ecor

Such a corollary, which is more significant if $\tau>0$,
as we have assumed,
can be applied to the problem
obtained by \lineariz ing \Pbl\ around its solution.
Therefore, it is \pier{useful} in the control problem associated to \Pbl\
we are going to discuss in a forthcoming paper.
Our second regularity result that \pier{deals} with the general case is the following

\Bthm
\label{Bddness}
In addition to the assumptions of Theorem~\ref{Regularity},
suppose that $\tau>0$ and that
\Beq
  g \in \LQ\infty , \quad
  \gG \in \LS\infty
  \aand
  \betaz(\yz) \in \LQ\infty .
  \label{hpbdd}
\Eeq
Then, the solution to problem \pbl\ also satisfies
\Beq
  w \in \L\infty\Hdue \subset \LQ\infty 
  \aand
  \xi \in \LQ\infty .
  \label{bddness}
\Eeq
\Ethm

The regularity result just stated has an interesting consequence
in the case of operators $\beta$ and~$\betaG$ satisfying the following assumptions
\Bsist
  && \hbox{$D(\beta)$ is an open interval~$I$}
  \aand
  D(\betaG) = D(\beta) .
  \label{hpstrutturater}
\Esist
The first \eqref{hpstrutturater} is fulfilled if $\Beta$ is, 
for instance, the everywhere smooth potential \eqref{regpot} 
or the logarithmic potential~\eqref{logpot}.
On the contrary, potentials whose convex part is an indicator function are excluded.
We observe that, if $I$ is not the whole of $\erre$ and $r_0$ is an end-point of it,
then $\betaz$ has an infinite limit at~$r_0$
since the interval $I$ is open.
Due to the second \pier{condition in}~\eqref{hpstrutturater},
the same remarks hold for $\betaGz$.
The result we state easily follows from~\eqref{hpstrutturater},
on account of \eqref{ybdd} (to~be used if $I$ is unbounded)
and the second \pier{property in}~\eqref{bddness}.
Therefore, we do not prove~it.

\Bcor
\label{StrongBddness}
In addition to the assumptions of Theorem~\ref{Regularity},
suppose that $\tau>0$ and that \eqref{hpbdd} and \eqref{hpstrutturater} are satisfied.
Moreover, assume that $\lam=1$ and $\lamG=1$.
Then, the following conclusions hold true:
$i)$~for the solution $(y,\yG,w, \xi,\xiG)$ to problem \Pbl\ we~have 
\Beq
  y(x,t) \in K 
  \quad \hbox{\aaQ\ and some compact subset $K\subset I$}.
  \label{compactvualed}
\Eeq
In particular, even $\xiG$ is bounded.
$ii)$ Assume that $\beta$ and $\betaG$ are single-valued $C^1$ functions.
Then, the functions $\beta'(y)$ and $\betaG'(\yG)$ are bounded as well.
$iii)$~Assume that $\beta$, $\betaG$, $\pi$ and $\piG$ are of class~$C^2$, in addition.
Then
\Beq
  \beta'(y) + \pi'(y) \in \L\infty V \cap \LQ\infty
  \aand
  \betaG'(y) + \piG'(y) \in \L\infty\VG \cap \LS\infty .
  \non
\Eeq
\Ecor

The rest of the section is devoted to recall some facts that are well known
and to introduce some notation that is widely used in the sequel.
First of all, we often owe to the Young inequality
(mainly with $p=p'=2$, thus with $\delta^{-p'/p}=\delta^{-1}$)
\Beq
  ab \leq \frac \delta p \, a^p + \frac {\delta^{-p'/p}} {p'} \, b^{p'}
  \quad \hbox{for every $a,b\geq 0$, $\delta>0$ and $p>1$}
  \label{young}
\Eeq
where $p':=p/(p-1)$,
and to the \Holder\ inequality.
Moreover, we account for the well-known embeddings and the related inequalities,
as well as the Poincar\'e inequality, namely
\begin{align}
  & \norma v_\infty \leq C \norma v_{\Hx2}
  \quad \hbox{for every $v\in\Hx2$}
  \label{sobolevHdue}
  \\
  & \norma v_\infty \leq C \norma v_{\HxG2}
  \quad \hbox{for every $v\in\HxG2$}
  \label{sobolevHdueG}
  \\
  & \normaV v^2 \leq C \bigl( \normaH{\nabla v}^2 + |v_\Omega|^2 \bigr)
  \quad \hbox{for every $v\in V$}
  \label{poincare}
\end{align}
where $C$~depends \pier{only on~$\Omega$}.
Furthermore, we observe that the identity 
$\normaH v^2=\<v,v>$ for $v\in V$ easily implies the inequality
\Beq
  \normaH v^2
  \leq \delta \normaH{\nabla v}^2
  + \cd \normaVp v^2
  \quad \hbox{for every $v\in V$}
  \label{interpolaz}
\Eeq
for every $\delta>0$ and some constant $\cd$ depending on $\delta$ and $\Omega$ as well.
Next, we recall a tool that is generally used 
in the context of problems related to the Cahn\pier{--}Hilliard equations.
We define
\Beq
  \dom\calN := \graffe{\vstar\in\Vp: \ \vstar_\Omega = 0}
  \aand
  \calN : \dom\calN \to \graffe{v \in V : \ v_\Omega = 0}
  \label{predefN}
\Eeq
by setting for $\vstar\in\dom\calN$
\Beq
  {\calN\vstar \in V, \quad
  (\calN\vstar)_\Omega = 0 ,
  \aand
  \iO \nabla\calN\vstar \cdot \nabla z = \< \vstar , z >
  \quad \hbox{for every $z\in V$}}
  \label{defN}
\Eeq
i.e., $\calN\vstar$ is the solution $v$ to the \generaliz ed Neumann problem for $-\Delta$
with datum~$\vstar$ that satisfies~$\pier{\vstar_{\,\Omega}=0}$.
As $\Omega$ is bounded, smooth, and connected,
it turns out that \eqref{defN} yields a well-defined isomorphism,
which satisfies
\Beq
  \< \ustar , \calN \vstar >
  = \< \vstar , \calN \ustar >
  = \iO (\nabla\calN\ustar) \cdot (\nabla\calN\vstar)
  \quad \hbox{for $\ustar,\vstar\in\dom\calN$}.
  \label{simmN}
\Eeq
Moreover, if we define $\normaVp\cpto:\Vp\to[0,+\infty)$ by the formula
\Beq
  \normaVp\vstar^2
  := \normaH{\nabla\calN(\vstar-(\vstar)_\Omega)}^2
  + |(\vstar)_\Omega)|^2
  \quad \hbox{for $\vstar\in\Vp$}
  \label{defnormaVp}
\Eeq
it is \sfw\ to prove that $\normaVp\cpto$ is a norm
that makes $\Vp$ a Hilbert space.
It follows that $\normaVp\cpto$ is equivalent to the usual dual norm
by the open mapping theorem
and it thus can be used as a norm in~$\Vp$.
It follows that
\Beq
  | \< \vstar , v >|
  \leq C \normaVp\vstar \normaV v
  \quad \hbox{for every $\vstar\in\Vp$ and $v\in V$}
  \non
\Eeq
where $C$ depends only on~$\Omega$.
Note that
\Beq
  \< \vstar , \calN\vstar >
  = \normaVp\vstar^2
  \quad \hbox{for every $\vstar\in\dom\calN$}
  \label{defnormavp}
\Eeq
by \accorpa{simmN}{defnormaVp}.
Finally, owing to \eqref{simmN} once more, we see~that
\Beq
  2 \< \dt\vstar(t) , \calN\vstar(t) >
  = \frac d{dt} \iO |\nabla\calN\vstar(t)|^2
  = \frac d{dt} \, \normaVp{\vstar(t)}^2
  \quad \aat
  \label{dtcalN}
\Eeq
for every $\vstar\in\H1\Vp$ satisfying $\pier{\vstar(t)_\Omega} =0$ 
for every $t\in[0,T]$.

We conclude this section by stating a general rule
we use as far as constants are concerned,
in order to avoid a boring notation.
Throughout the paper,
the small-case symbol $c$ stands for different constants which depend only
on~$\Omega$, on the final time~$T$, and on the constants and the norms of
the functions involved in the assumptions of our statements.
In particular, $c$~is independent of the approximation parameter~$\eps$
we introduce later~on.
A~notation like~$\cd$ allows the constant to depend on the positive
parameter~$\delta$, in addition.  
Hence, the meaning of $c$ and $\cd$ might
change from line to line and even in the same chain of inequalities.  
On the contrary, we use capital letters 
to~denote precise constants which we could refer~to
(see, e.g.,~\eqref{sobolevHdue}).

%%%%%%%%%%%%%%%%%%%%%%%%%%%%%%%%%%%%%%%%%%%%%%%%%%%%%%%%%%%%%%%%%%%%%%%%

\section{Uniqueness and continuous dependence}
\label{UniqContDip}
\setcounter{equation}{0}

\pier{This section is devoted to the proof of Theorem~\ref{ContDip}.
We closely follow \cite[Thm.~1]{GiMiSchi}
and} just adapt the argument used there.
For convenience, we set $y:=y_1-y_2$ and similarly define
$\yG$, $w$, $\xi$, $\xiG$, $g$, $\gG$ and~$\yz$.
As the initial data have the same mean value,
by Remark~\ref{Conserved} applied to $y_i$ for $i=1,2$,
we see that $y(t)$ has zero mean value 
and thus belongs to the domain of $\calN$ for every~$t\in[0,T]$.
Therefore, we can write \eqref{prima} at any time $s$ for both solutions
and test the difference by~$\calN y(s)$.
Then, we integrate over~$(0,t)$ with respect to~$s$,
where $t\in(0,T]$ is arbitrary.
At the same time, we write \eqref{seconda} for both solutions
and \pier{take~$-y$ as test function}.
Finally, we add the obtained equalities to each other.
We~have
\Bsist
  && \iot \< \dt y(s) , \calN y(s) > \, ds
  + \intQt \nabla w \cdot \nabla\calN y
  - \intQt wy
  \non
  \\
  && \quad {}
  + \frac \tau 2 \iO |y(t)|^2
  - \frac \tau 2 \iO |\yz|^2
  + \frac 12 \iG |\yG(t)|^2
  - \frac 12 \iG |\yz\suG|^2
  \non
  \\
  && \quad {}
  + \intQt |\nabla y|^2
  + \intSt |\nablaG\yG|^2
  + \intQt \xi y
  + \intSt \xiG \yG
  \non
  \\
  && = \intQt \bigl( \lam(\pi(y_2) - \pi(y_1)) + g \bigr) \, y
   + \intSt \bigl( \lamG(\piG(y_{\Gamma\!,2})) - \piG(y_{\Gamma\!,1}) + \gG\bigr) \, \yG \,.
  \non
\Esist
Now, we transform the first term on the \lhs\ with the help of \eqref{dtcalN}
and cancel the next two integrals \pier{thanks} to~\eqref{defN}.
Moreover, we neglect the last two integrals on the \lhs\
since they are nonnegative for $\beta$ and $\betaG$ are monotone.
Finally, we \pier{exploit} assumptions \accorpa{hppilip}{hplam}
and \pier{use} the elementary Young inequality.
We~obtain
\Bsist
  && \frac 12 \, \normaVp{y(t)}^2
  - \frac 12 \, \normaVp\yz^2
  + \frac \tau 2 \, \normaH{y(t)}^2
  - \frac \tau 2 \, \normaH\yz^2
  + \frac 12 \, \normaHG{\yG(t)}^2
  - \frac 12 \, \normaHG{\yz\suG}^2
  \non
  \\
  && \quad {}
  + \intQt |\nabla y|^2
  + \intSt |\nablaG\yG|^2
  \non
  \\
  && \leq c \intQt |y|^2
  + \frac 14 \intQt |g|^2
  + c \intSt |\yG|^2
  + \frac 14 \intSt |\gG|^2 .
  \non
\Esist
At this point, we \pier{take advantage of \eqref{interpolaz} to infer that}
\Beq
  \intQt |y|^2
  \leq \delta \intQt |\nabla y|^2
  + \cd \iot \normaVp{y(s)}^2 \, ds.
  \non
\Eeq
Therefore, it \pier{suffices} to choose \pier{$\delta <1 $} 
and apply the Gronwall lemma to obtain~\eqref{contdip}.
The sentence of the statement regarding partial uniqueness easily follows,
as we show at once.
Clearly, \eqref{contdip} with the same data implies 
$y=0$ and $\yG=0$ with the notation we have introduced at the beginning, 
so that the difference of the equation~\eqref{seconda} 
simply reduces~to
\Beq
  \iO w(t) \, v
  = \iO \xi(t) \, v
  + \iG \xiG(t) \, v
  \quad \hbox{\aat\ and every $v\in\calV$} .
  \label{quasiunicita}
\Eeq
\pier{In view of \eqref{prima}, it turns out that $w$ is a function depending only on time (i.e., \aat\ $w(t)$ is constant in $\Omega$). Taking now test functions $v
\in {\cal D} (\Omega)$ in \eqref{quasiunicita}, we easily infer that $\xi (t) = w(t) $ \aat. Next, letting $v$ vary in $\calV $ we also deduce that $\xiG(t)=0$ \aat. } Assume now $\beta$ to be single-valued. Then, we have $\xi=0$ in~\eqref{quasiunicita} \pier{as  well}.
We immediately \pier{conclude} that $w=0$ and the proof is complete.

%%%%%%%%%%%%%%%%%%%%%%%%%%%%%%%%%%%%%%%%%%%%%%%%%%%%%%%%%%%%%%%%%%%%%%%%

\section{Existence and regularity}
\label{ExistenceRegularity}
\setcounter{equation}{0}

In this section, we prove our existence and regularity results.
\pier{The method we use} relies on a \regulariz ation depending on the parameter $\eps\in(0,1)$ that will tend to zero.
Namely, we introduce the approximating problem
of finding a pair $(\yeps,\weps)$ such that 
a suitably corresponding quintuplet $(\yeps,\yGeps,\weps,\xieps,\xiGeps)$
solves the system obtained by replacing the graphs $\beta$ and~$\betaG$
by the everywhere defined functions $\betaeps$ and~$\betaGeps$ we make precise below
and $g$ by a suitably \regulariz ed datum.
For clarity, we write both the construction of the corresponding quintuplet
and the \regulariz ed problem, at once:
\begin{align}
  & \yGeps(t) := \yeps(t)\suG \quad \aat, \quad
  \xieps := \betaeps(\yeps)
  \aand
  \xiGeps := \betaGeps(\yGeps)
  \label{relazionieps}
  \\
  \noalign{\smallskip}
  & \< \dt\yeps(t) , v >
  + \iO \nabla\weps(t) \cdot \nabla v = 0
  \label{primaeps}
  \\
  \noalign{\smallskip}
  & \iO \weps(t) v
  = \taueps \iO \dt\yeps(t) \, v
  + \iG \dt\yGeps(t) \, v
  + \iO \nabla\yeps(t) \cdot \nabla v
  + \iG \nablaG\yGeps(t) \cdot \nablaG v
  \qquad
  \non
  \\
  & \quad {}
  + \iO \bigl( \xieps(t) + \lam(t) \, \pi(\yeps(t)) - \geps(t) \bigr) \, v
  + \iG \bigl( \xiGeps(t) + \lamG(t) \, \piG(\yGeps(t)) - \gG(t) \bigr) \, v
  \label{secondaeps}
  \\
  \noalign{\smallskip}
  & \yeps(0) = \yz 
  \label{cauchyeps}
\end{align}
\Accorpa\Pbleps primaeps cauchyeps
where \eqref{primaeps} and \eqref{secondaeps} are required to hold
for every $v\in V$ and every $v\in\calV$, respectively,~and
\Beq
  \taueps :=\tau
  \quad \hbox{if $\tau>0$}
  \aand
  \taueps := \eps
  \quad \hbox{if $\tau=0$} .
  \label{deftaueps}
\Eeq

Thus, we first solve problem \Pbleps\ in \pier{the} proper functional framework.
Then, we perform a number of a~priori estimates and use compactness and monotonicity techniques
that ensure that the $\eps$-solution converges 
to a solution to the original problem in a proper topology as $\eps$ tends to zero.
Due to uniqueness, the whole family of approximating solution will converge,
even though it is necessary to take convergent subsequences, in principle.
The power of the estimates we can derive (thus, the topology of the convergence that follows)
depends on the assumptions on the data we can account for, i.e., on the theorem we want to prove.
We start with Theorem~\ref{Existence} and suppose that just \HPdati\ are fulfilled.
However, the whole argument partially works for the proof of the regularity results as well.
Just further a~priori estimates are necessary for the latter, indeed.
For the approximating solution we postulate the following regularity
\Bsist
  && \yeps \in \H1H \cap \L\infty V \cap \L2\Hdue
  \label{piuregyeps}
  \\
  && \yGeps \in \H1\HG \cap \L\infty\VG \cap \L2\HdueG
  \label{piuregyGeps}
  \\
  && \weps \in \L2V .
  \label{regweps}
\Esist
\Accorpa\Regsoluzeps piuregyeps regweps
Thus, we look for a pair $(\yeps,\weps)$
such that the quintuplet $(\yeps,\yGeps,\weps,\xieps,\xiGeps)$
defined by \eqref{relazionieps} satisfies \Regsoluzeps\ and solves \Pbleps.

Let us come to the definition of the \regulariz ed monotone operators.
Inspired by~\cite{CaCo}, for each graph,
we take the proper Yosida \regulariz ation (see~\eqref{yosida}), namely
\Beq
  \betaeps := \betaeps^Y
  \aand
  \betaGeps := (\betaG)_{\eta\eps}^Y
  \label{defbetaeps}
\Eeq
where $\eta$ is the same as in~\eqref{hpcompatib}.
Such a choice yields (see~\cite[Lemma 4.4]{CaCo})
\Beq
  |\betaeps(r)| \leq \eta |\betaGeps(r)| + C
  \quad \hbox{for every $r\in\erre$ and $\eps\in(0,1)$}
  \label{compatibeps}
\Eeq
where $C$ is some positive constant 
and $\eta$ still is the same as in~\eqref{hpcompatib}.
We also define for convenience
\Beq
  \Betaeps(r) := \int_0^r \betaeps(s) \, ds
  \aand
  \BetaGeps(r) := \int_0^r \betaGeps(s) \, ds
  \quad \hbox{for $r\in\erre$}
  \label{defBetaeps}
\Eeq
and recall that the set of properties
\Bsist
  \hskip-1cm && \hbox{$0\leq\Betaeps(r)\leq\Beta(r)$, \ $\Betaeps(r)\Neto\Beta(r)$ monotonically as $\eps\Seto0$}
  \label{propYosida}
  \\
  \hskip-1cm && \hbox{$|\betaeps(r)|\leq|\betaz(r)|$, \ $\betaeps(r)$ tends to $\betaz(r)$ monotonically as $\eps\Seto0$}
  \label{propyosida}
\Esist
and the analogue for $\BetaGeps$ and $\betaGeps$ hold true (see, e.g., \cite[Prop.~2.11, p.~39]{Brezis}).
Notice that inequality \eqref{compatibeps} implies
\Beq
  \Betaeps(r) \leq \eta \,\BetaGeps(r) + C |r|
  \quad \hbox{for every $r\in\erre$ and $\eps\in(0,1)$}
  \label{dacompatibeps}
\Eeq
since $\betaeps(0)=\betaGeps(0)=0$ by~\eqref{hpBeta}, whence $\betaeps$ and $\betaGeps$ have the same sign.
For the approximating datum~$\geps$,
we require the following regularity and convergence properties
\begin{align}
  & \geps \in \H1H
  \aand
  \geps \to g
  \quad \hbox{strongly in $\L2H$}
  \quad \hbox{if $\tau>0$}
  \label{reggeps}
  \\
  & \geps = g 
  \quad \hbox{if $\tau=0$} .
  \label{nogeps}
\end{align}
\Accorpa\Datieps reggeps nogeps
To start with our program, we first have to solve the approximating problem.

\Bthm
\label{Existenceeps}
Assume \pier{\accorpa{hpBeta}{defbeta}}, \HPdati, \eqref{deftaueps}, \eqref{defbetaeps} and \Datieps.
Then, there \pier{is} a unique pair $(\yeps,\weps)$ such that
the corresponding quintuplet $(\yeps,\yGeps,\weps,\xieps,\xiGeps)$ satisfies 
the regularity given by \Regsoluzeps\ and solves problem \Pbleps.
\Ethm

\Bdim
Uniqueness follows from Theorem~\ref{ContDip} as a particular case.
As far as existence is concerned,
we can first quote \cite[Thm.~3]{GiMiSchi},
even though $\pi$ and $\piG$ are \pier{multiplied} by coefficients
that depend on $x$ and~$t$.
Indeed, minor changes in the proof are sufficient 
to adapt the argument and obtain the \generaliz ed result we need here, namely,
the existence of a solution satisfying
\begin{multline*}
  \pier{\yeps \in  \H1H \cap \L\infty V  , \quad
  \yGeps \in \H1\HG \cap \L\infty\VG}  \\
    \aand
  \weps \in \L2V .
  \non
\end{multline*}
For completeness, we note that the cited result regards 
everywhere defined monotone operators $\beta$ and $\betaG$
satisfying suitable growth conditions 
(that include the sublinear case of the operators $\betaeps$ and~$\betaGeps$)
instead of compatibility conditions,
and data satisfying assumptions that are implied by \HPdati\ and \Datieps.
Thus, we only have to show that the solution we have found is smoother than expected, 
that is, it satisfies \accorpa{piuregyeps}{piuregyGeps}.
In this direction, we could quote \cite[Rem.~12]{GiMiSchi}\pier{; however,}
we detail the argument for the reader's convenience.
We observe that the variational equation \eqref{secondaeps} implies
that $\yeps$ solves the partial differential equation
\Beq
  \weps
  = \taueps \dt\yeps - \Delta\yeps + \betaeps(\yeps) + \lam \, \pi(\yeps) - \geps
  \quad \hbox{in $Q$}
  \label{pdeeps}
\Eeq
at least in the sense of distributions.
Due to the regularity assumed for $\geps$
and the \Lip\ continuity of~$\betaeps$,
all the terms of \eqref{pdeeps} but $\Delta\yeps$ belong to~$\L2H$.
By comparison, we deduce that $\Delta\yeps\in\L2H$.
On the other hand, $\yGeps\in\L2{\pier{\VG}}$.
Thus, the elliptic regularity theory
yields
$\yeps\in\L2{\Hx{3/2}}$, whence also $\dn\yeps\suG\in\L2\HG$ 
(see, e.g., \cite[Thms.~7.4 and~7.3, pp.\ 187-188]{LioMag}
or \cite[Thm.~3.2, p.~1.79, and Thm.~2.27, p.~1.64]{BreGil}).
In particular, all the terms of the integration by part formula for the Laplace operator
are functions and we deduce that the variational equation \eqref{secondaeps} also implies
\Beq
  \dn\yeps\suG + \dt\yGeps - \DeltaG\yGeps + \betaGeps(\yGeps) + \lamG \piG(\yGeps) = \gG
  \quad \hbox{on $\Sigma$}
  \label{bceps}
\Eeq
at least in a \generaliz ed sense, in principle.
Arguing as before, 
we see that $\Delta\yGeps\in\L2\HG$, whence
$\yGeps\in\L2\HdueG$
by~the boundary version of the elliptic regularity theory. 
Coming back to~$\yeps$, we infer that $\yeps\in\L2\Hdue$.
\Edim

At this point, we start estimating.
We are widely inspired by the techniques used in~\cite{GiMiSchi}.
However, since our assumptions and statements are different,
it is necessary to detail the argument.
We remind the reader that our assumptions on the data reduce to \HPdati\ and \Datieps.
Moreover we recall the definition \eqref{hpinterno} of~$\mz$
and observe that Remark~\ref{Conserved} obviously applies to the approximating problem as well,~i.e.
\Beq
  (\dt\yeps(t))_\Omega = 0
  \quad \aat
  \aand
  \yeps(t)_\Omega = \mz
  \quad \hbox{for every $t\in[0,T]$} .
  \label{conservedeps}
\Eeq
Furthermore, we recall the \pier{useful} inequalities
\begin{align}
  \betaeps(r)(r-\mz) \geq \deltaz |\betaeps(r)| - \Cz
  \aand
  \betaGeps(r)(r-\mz) \geq \deltaz |\betaGeps(r)| - \Cz
  \qquad 
  \non
  \\
  \quad \hbox{for some $\deltaz,\Cz>0$ and every $r\in\erre$ and $\eps\in(0,1)$}
  \label{mirzelik}
\end{align}
which hold whenever $\Beta(0)=\BetaG(0)=0$ 
and $\mz$ lies in the interior of the domains of $\beta$ and~$\betaG$
(and $\deltaz$ and $\Cz$ depend on the position of~$\mz$),
thus under our assumptions
(see \eqref{hpBeta}, \eqref{hpinterno} and the inclusion in~\eqref{hpcompatib}).
For a detailed proof, see \cite[p.~908]{GiMiSchi}; 
see also \cite[Appendix, Prop.~A.1]{MirZelik}.

\step
First a priori estimate

We write \eqref{primaeps} at the time $s$
and take $v=\calN(\yeps(s)-\mz)$, which is meaningful by~\eqref{conservedeps}
(see~\eqref{predefN}).
Then, we integrate over $(0,t)$ with respect to~$s$,
where $t$ is arbitrary in~$(0,T]$.
At the same time, we analogously behave with \eqref{secondaeps}
by choosing $-(\yeps-\mz)$ as a test function.
Then, we sum the obtained equalities to each other
and use \eqref{dtcalN} and \eqref{defN}
in order to transform the first integral we get and to cancel the next two ones.
Finally, we add \pier{two additional} terms to both sides for convenience.
We obtain
\begin{align}
  & \frac 12 \, \normaVp{\yeps(t)-\mz}^2
  + \frac \taueps 2 \iO |\yeps(t)-\mz|^2
  + \frac 12 \iG |\yGeps(t)-\mz|^2
  \non
  \\
  & \quad {}
  + \intQt |\nabla\yeps|^2
  + \intSt |\nablaG\yGeps|^2
  \non
  \\
  & \quad {}
  + \intQt \bigl( \betaeps(\yeps) - \betaeps(\mz) \bigr) (\yeps-\mz)
  + \intSt \bigl( \betaGeps(\yGeps(t)) - \betaGeps(\mz) \bigr) (\yGeps-\mz)
  \non
  \\
  & = \frac 12 \, \normaVp{\yz-\mz}^2
  + \frac \taueps 2 \iO |\yz-\mz|^2
  + \frac 12 \iG |\yz\suG-\mz|^2
  - \betaGeps(\mz) \intSt (\yGeps-\mz)
  \non
  \\
  & \quad {}
  + \intQt \bigl( \geps - \lam \, \pi(\yeps) \bigr) (\yeps-\mz)
  + \intSt \bigl( \gG - \lamG \piG(\yGeps) \bigr) (\yeps-\mz) .
  \label{perprimastima}
\end{align}
The last two integrals on the \lhs\ are nonnegative by monotonicity
and one term on the \rhs\ can be dealt with this~way
\Beq
  - \betaGeps(\mz) \intSt (\yGeps-\mz)
  \leq c \, |\betaGeps(\mz)|^2 + \intSt |\yGeps-\mz|^2
  \leq c + \intSt |\yGeps-\mz|^2
  \non
\Eeq
thanks to the analogue of \eqref{propyosida} for $\betaGeps$ and to assumption \eqref{hpinterno} on~$\mz$.
On the other hand, we can \pier{take \accorpa{hppilip}{hplam} and \eqref{hpg} into account} and apply \eqref{interpolaz} to~$\yeps-\mz$, 
in order to estimate the sum of the last two terms on the \rhs\ of \eqref{perprimastima} as follows
\Bsist
  &&\pier{\intQt \bigl( \geps - \lam \, \pi(\yeps) \bigr) (\yeps-\mz)
  + \intSt \bigl( \gG - \lamG \piG(\yGeps) \bigr) (\yeps-\mz)}
  \non
  \\
  && \leq c
  + c \intQt |\yeps-\mz|^2 
  + c \intSt |\yGeps-\mz|^2 
  \non
  \\
  && \leq c
  + \frac 12 \intQt |\nabla\yeps|^2 
  + c \iot \normaVp{\yeps(s)-\mz}^2 \, ds
  + c \intSt |\yGeps-\mz|^2 .
  \non
\Esist
At this point, we combine \eqref{perprimastima} with the above inequalities, 
rearrange a little and apply the Gronwall lemma.
Then, we easily eliminate $\mz$ in the estimate we obtain.
By~using \eqref{interpolaz} once more,
we recover the $L^2$ norm of~$\yeps$ through the norm of its gradient
and eventually conclude~that
\Beq
  \norma\yeps_{\L\infty\Vp\cap\L2V}
  + \taueps^{1/2} \norma\yeps_{\L\infty H}
  + \norma\yGeps_{\L\infty\HG\cap\L2\VG}
  \leq c \,.
  \label{primastima}
\Eeq

\step
Consequence

In view of~\accorpa{hppilip}{hplam}, we immediately deduce~that
\Beq
  \norma{\lam\,\pi(\yeps)}_{\L2H}
  + \norma{\lamG\piG(\yGeps)}_{\L2\HG}
  \leq c \,.
  \label{stimapipiG}
\Eeq
Moreover, by \pier{virtue of}~\eqref{hpdxlam}\ if $\tau=0$, we have
\Beq
  |\nabla(\lam\,\pi(\yeps))|
  = |\pier{\pi(\yeps)} \, \nabla\lam + \lam \, \pi'(\yeps) \nabla\yeps|
  \leq \pier{c\left( |\nabla\lam| \, |\yeps| +  |\nabla\yeps| \right)} \,.
  \non
\Eeq
Hence, by also accounting for the \Holder\ inequality 
and the continuous embedding $V\subset\Lx 6$, we deduce~that
\Beq
  \norma{\nabla(\lam\,\pi(\yeps))}_{\L2H}
  \leq c \norma{\nabla\lam}_{\L\infty{\Lx3}} \norma\yeps_{\L2{\Lx6}}
  + c \norma\yeps_{\L2V}
  \leq c 
  \non
\Eeq
and conclude that
\Beq
  \norma{\lam\,\pi(\yeps)}_{\L2V} \leq c 
  \quad \hbox{if $\tau=0$}.
  \label{stimapiV}
\Eeq

\step
Second a priori estimate

We recall~\eqref{conservedeps} and test \eqref{primaeps} and \eqref{secondaeps} 
written at the time~$s$ by~$\calN(\dt\yeps(s))$ and $-\dt\yeps(s)$, respectively.
We \pier{sum the obtained equalities} and integrate over~$(0,t)$.
Then, \pier{recalling} \eqref{defN} and \eqref{defnormavp} once more\pier{, we have that}
\Bsist
  && \iot \normaVp{\dt\yeps(s)}^2 \, ds
  + \taueps \intQt |\dt\yeps|^2
  + \intSt |\dt\yGeps|^2
  \non
  \\
  && \quad {}
  + \frac 12 \iO |\nabla\yeps(t)|^2
  + \frac 12 \iG |\nablaG\yGeps(t)|^2
  + \iO \pier{\Betaeps} (\yeps(t))
  + \iG \pier{\BetaGeps} (\yGeps(t))
  \non
  \\
  && = %\frac 12 \, \normaVp\yz^2
   \frac 12 \iO |\nabla \yz|^2
  + \frac 12 \iG |\nablaG\yz\suG|^2
  + \iO \pier{\Betaeps} (\yz)
  + \iG \pier{\BetaGeps} (\yz\suG)
  \non
  \\
  && \quad {}
  + \intQt \bigl( \geps - \lam \, \pi(\yeps) \bigr) \dt\yeps
  + \intSt \bigl( \gG - \lamG \piG(\yGeps) \bigr) \dt\yGeps.
  \label{persecondastima}
\Esist
Note that all the terms on the \lhs\ are nonnegative 
(cf.~\pier{\eqref{propYosida}}) and recall that \eqref{hpyz} holds.
\pier{Then, the upper inequalities in \eqref{propYosida}, holding for $\BetaGeps$ and $\betaG$ as well, allow us to infer
$$
\iO \pier{\Betaeps} (\yz) + \iG \pier{\BetaGeps} (\yz\suG) \leq 
\iO \pier{\Beta} (\yz) + \iG \pier{\BetaG} (\yz\suG) \leq c. 
$$
}%
Thus, just the last two integrals on the \rhs\ need some treatment.
By \eqref{hpg}, \eqref{stimapipiG} and the Young inequality~\eqref{young},
we immediately have
\Beq
  \intSt \bigl( \gG - \lamG \piG(\yGeps) \bigr) \dt\yGeps
  \leq c
  + \frac 12 \, \norma{\dt\yGeps}_{\L2\HG}^2 .
  \non
\Eeq
On the contrary, the analogous integral over $Q_t$ is more delicate
and we distinguish the cases $\tau>0$ and $\tau=0$.
In the former, we~write
\Beq
  \intQt \bigl( \geps - \lam\,\pi(\yeps) \bigr) \dt\yeps
  \leq \norma{\geps-\lam\,\pi(\yeps)}_{\L2H}^2
  + \frac \tau 2 \, \norma{\dt\yeps}_{\L2H}^2 
  \leq c
  + \frac \tau 2 \, \norma{\dt\yeps}_{\L2H}^2 .
  \non
\Eeq
Hence, the last term can be absorbed by the \lhs\ of \eqref{persecondastima}.
If $\tau=0$, we have $\geps=g$ and can account for~\eqref{hpgtauzero}.
\pier{In view of}~\eqref{primastima}, \eqref{stimapiV}
and the interpolation inequality~\eqref{interpolaz},
we obtain
\Bsist
  && \intQt \geps \dt\yeps
  = \iO g(t) \yeps(t)
  - \iO g(0) \yz
  - \intQt \dt g \, \yeps
  \non
  \\
  && \leq \norma g_{\L\infty H}^2 + \normaH{\yeps(t)}^2
  + c
  + \pier{\norma {\dt g}_{\L2H}}^2 + \norma\yeps_{\L2H}^2
  \non
  \\
  && \leq \frac 14 \iO |\nabla\yeps(t)|^2 
  + c \normaVp{\yeps(t)}^2 + c
  \leq \frac 14 \iO |\nabla\yeps(t)|^2 
  + c 
  \non
\Esist
\pier{and}
\Beq
  - \intQt \lam \, \pi(\yeps) \dt\yeps
  \leq \frac 12 \iot \normaVp{\dt\yeps(s)}^2 \, ds
  + c \norma{\lam\,\pi(\yeps)}_{\L2V}^2
  \leq \frac 12 \iot \normaVp{\dt\yeps(s)}^2 \, ds
  +c \,.
  \non
\Eeq
Thus, the \lhs\ of \eqref{persecondastima} dominates also in this case.
Therefore, we conclude that
\begin{align}
  & \norma\yeps_{\H1\Vp\cap\L\infty V}
  + \taueps^{1/2} \norma\yeps_{\H1H}
  + \norma\yGeps_{\H1\HG\cap\L\infty\VG}
  \non
  \\
  & \quad {}
  + \norma{\Betaeps(\yeps)}_{\L\infty\Luno}
  + \norma{\BetaGeps(\yGeps)}_{\L\infty\LunoG}
  \leq c \,.
  \label{secondastima}
\end{align}

\step
Third a priori estimate

From \eqref{secondastima}, we first deduce an estimate for $\nabla\weps$ 
by using~\eqref{primaeps} with $v=\weps(t)-(\weps(t))_\Omega$
and the Poincar\'e inequality~\eqref{poincare}.
We have \aat
\Bsist
  && \iO |\nabla\weps(t)|^2 
  = \iO |\nabla \bigl( \weps(t)-(\weps(t))_\Omega \bigr)|^2 
  = \< \dt\yeps(t) , \weps(t)-(\weps(t))_\Omega >
  \non
  \\
  && \leq {}\pier{C}\, \normaVp{\dt\yeps(t)} \, 
  \normaV{\weps(t)-(\weps(t))_\Omega}
  \leq \frac 12 \iO |\nabla\weps(t)|^2 
  + c \, \normaVp{\dt\yeps(t)}^2 \,.
  \non
\Esist
Hence, \eqref{secondastima} immediately yields
\Beq
  \norma{\nabla\weps}_{\L2H} \leq c \,.
  \label{stimanablaw}
\Eeq
In order to recover the full $V$-norm, we have to estimate the mean value.
Inspired by \cite[p.~908]{GiMiSchi}, we test equations \eqref{primaeps} and \eqref{secondaeps}
as we did for our first estimate, i.e., by $\calN(\yeps-\mz)$ and $-(\yeps-\mz)$, respectively,
but we do not integrate with respect to time.
Also in the present case, two terms cancel thanks to~\eqref{defN}.
Thus, \aat\ (and we avoid writing $t$ in the sequel, for brevity) we~have
\Bsist
  && \iO |\nabla\yeps|^2
  + \iG |\nablaG\yGeps|^2
  + \iO \xieps (\yeps-\mz)
  + \iG \xiGeps (\yGeps-\mz)
  \non
  \\
  && = F_\eps
  := - \< \dt\yeps , \calN(\yeps-\mz) >
  + \iO \bigl( \geps - \lam \, \pi(\yeps) - \taueps \dt\yeps \bigr) (\yeps-\mz)
  \non
  \\
  && \quad {}
  + \iG \bigl( \gG - \lamG \piG(\yGeps) - \dt\yGeps \bigr) (\yGeps-\mz) .
  \label{perstimamediaw}
\Esist
Now, we account for \eqref{mirzelik} and deduce that
\Beq
  \iO \xieps (\yeps-\mz)
  + \iG \xiGeps (\yGeps-\mz)
  \geq \deltaz \iO |\betaeps(\yeps)|
  + \deltaz \iG |\betaGeps(\yGeps)|
  - c \,.
  \label{damirzelik}
\Eeq
On the other hand, recalling \eqref{defnormaVp} and that $\yeps$ and $\yGeps$
are bounded in $\L\infty V$ and in $\L\infty\VG$, respectively (see~\eqref{secondastima}),
we deduce \aat
\Bsist  
  && |F_\eps|
  \leq c \normaVp{\dt\yeps} \, \normaVp{\yeps-\mz}
  + \bigl( \normaH\geps + \normaH{\lam\,\pi(\yeps)} + \taueps \normaH{\dt\yeps} \bigr) \normaH{\yeps-\mz}
  \non
  \\
  && \qquad\quad {}
  + \bigl( \normaHG\gG + \normaHG{\lamG\piG(\yGeps)} + \normaHG{\dt\yGeps} \bigr) \normaHG{\yGeps-\mz} 
  \non
  \\
  && \leq c \normaVp{\dt\yeps}  
  + c \bigl( \normaH\geps + \normaH{\lam\,\pi(\yeps)} + \taueps \normaH{\dt\yeps} \bigr) 
  \non
  \\
  && \quad {}
  + c \bigl( \normaHG\gG + \normaHG{\lamG\piG(\yGeps)} + \normaHG{\dt\yGeps} \bigr) .
  \non
\Esist
Hence, $F_\eps$ is bounded in $L^2(0,T)$
thanks to \eqref{secondastima}, \eqref{stimapipiG} and assumptions \eqref{hpg}, \eqref{reggeps} on the data.
Therefore, even the integrals on the \rhs\ of \eqref{damirzelik} are estimated in~$L^2(0,T)$. 
By choosing $v=1$ in~\eqref{secondaeps}, we immediately deduce that
the same holds for the space integral of~$\weps$, \pier{whence
\Beq
  \bigl\|(\weps)_\Omega \bigr\|_{L^2(0,T)}
  \leq c 
  \label{stimamedia}
\Eeq
i.e., the mean value of $\weps$} is estimated in~$L^2(0,T)$.
By combining this and~\eqref{stimanablaw} and using the Poincar\'e inequality~\eqref{poincare}, 
we eventually infer~that
\Beq
  \norma\weps_{\L2V} \leq c \,.
  \label{stimaw}
\Eeq

\step 
Fourth a priori estimate

Our aim is to find \pier{a bound for~$\xieps$ in $L^2(Q)$}.
To this end, we \pier{take $v=\xieps$ in \eqref{secondaeps}} 
and integrate \pier{over $\Omega$. We have
\begin{align}
  & 
  \iO \betaeps'(\yeps) |\nabla\yeps|^2
  + \iG \betaeps'(\yGeps) |\nablaG\yGeps|^2
  + \iO |\xieps|^2
  + \iG \betaGeps(\yGeps) \betaeps(\yGeps)
  \non
  \\
  &= \iO \bigl( g - \lam\,\pi(\yeps) + \weps - \tau \dt\yeps \bigr) \xieps
  + \iG \bigl( \gG - \lamG\piG(\yGeps) - \dt\yGeps\bigr) \betaeps(\yGeps) .
  \label{perquartastima}
\end{align}
The first three terms on the \lhs\ are nonnegative.
For the last one, we recall the compatibility condition~\eqref{compatibeps}:
note that} the functions $\betaeps$ and $\betaGeps$ have the same sign
since they are \pier{non-decreasing and null} at~$0$ (due to~\eqref{hpBeta})\pier{. Then}, we \pier{deduce that
\Beq
  \iG \betaGeps(\yGeps) \betaeps(\yGeps)
  \geq \frac 1\eta \iG \bigl( |\betaeps(\yGeps)|^2 
  - C |\betaeps(\yGeps)| \bigr)
  \geq \frac 1 {2\eta} \iG |\betaeps(\yGeps)|^2 - c \,.
  \non
\Eeq
}%
Let us come to the \rhs\ of \eqref{perquartastima}. 
The first two terms are bounded
thanks to~\eqref{propYosida}, \eqref{dacompatibeps} and~\eqref{hpyz}.
Furthermore, \pier{using the Young inequality \eqref{young}}
we can estimate the last integrals as follows\pier{
\begin{align}
  & \iO \bigl( g - \lam\,\pi(\yeps) + \weps  - \tau \dt\yeps\bigr) \xieps
  \leq \frac 12 \iO |\xieps|^2 + \frac 12 \norma{g - \lam\,\pi(\yeps) + \weps  - \tau \dt\yeps}_H^2
  \non
  \\
  & \iG \bigl( \gG - \lamG\piG(\yGeps) - \dt\yGeps\bigr) \betaeps(\yGeps)
  \leq \frac 1 {4\eta} \iG |\betaeps(\yGeps)|^2 + \eta \, \norma{  \gG - \lamG\piG(\yGeps) - \dt\yGeps}_{\HG}^2
  \non
\end{align}
and remark that, thanks to \eqref{stimapipiG}, \eqref{secondastima} and \eqref{stimaw}, the last terms in the above inequalities are uniformly bounded in $L^1(0,T).$ Hence, by combining and integrating over $(0,T)$, we find out that}
\Beq
  \norma\xieps_{\L2H} \leq c .
  \label{quartastima}
\Eeq

\step
Consequence

By partially repeating the argument used in the proof of Theorem~\ref{Existenceeps}
and noting that each deduction has a corresponding estimate,
we derive the following chain of bounds
\Beq
  \norma{\Delta\yeps}_{\L2H} \leq c, \quad
  \norma\yeps_{\L2{\Hx{3/2}}} \leq c , \quad
  \norma{\dn\yeps\suG}_{\L2\HG} \leq c .
  \non
\Eeq
By comparison in \eqref{bceps} with the help of~\eqref{secondastima}, 
we conclude that
\Beq
  \norma{-\DeltaG\yGeps + \betaGeps(\yGeps)}_{\L2\HG} \leq c \,.
  \label{perstimabetaG}
\Eeq

\step
Fifth a priori estimate

By \eqref{perstimabetaG}, we can simply write 
\Beq
  -\DeltaG\yGeps + \betaGeps(\yGeps) = F_\eps \, , \ 
  \pier{\hbox{ with }} \ 
  \norma{F_\eps}_{\L2\HG} \leq c 
  \non
\Eeq
and multiply such an equation by~$\betaGeps(\yGeps)$.
We immediately obtain
\Beq
  \intS \betaGeps'(\yGeps) |\nablaG\yGeps|^2
  + \intS |\betaGeps(\yGeps)|^2
  = \intS F_\eps \, \betaGeps(\yGeps)
  \leq \frac 12 \intS |\betaGeps(\yGeps)|^2
  + c 
  \non
\Eeq
and infer that
\Beq
  \norma{\betaGeps(\yGeps)}_{\L2\HG} \leq c \,.
  \label{stimabetaGeps}
\Eeq

\step
Consequence

By entering the proof of Theorem~\ref{Existenceeps} once more
and arguing as above,
we deduce the following chain of bounds
\Beq
  \norma{\DeltaG\yGeps}_{\L2\HG} \leq c, \quad
  \norma\yGeps_{\L2\HdueG} \leq c , \quad
  \norma\yeps_{\L2\Hdue} \leq c \,.
  \non
\Eeq
Therefore, we conclude that
\Beq
  \norma\yGeps_{\L2\HdueG} \leq c 
  \aand
  \norma\yeps_{\L2\Hdue} \leq c 
  \label{stimeHdue}
\Eeq
whence also
\Beq
  \norma\yGeps_{\L2\LinftyG} \leq c 
  \aand
  \norma\yeps_{\L2\Linfty} \leq c 
  \label{dastimeHdue}
\Eeq
thanks to the continuous embeddings \accorpa{sobolevHdue}{sobolevHdueG}.

\step
Conclusion of the proof of Theorem \ref{Existence}

Recalling all the estimates we have obtained,
we see that the following convergence holds true
\Bsist
  \hskip-0.8cm & \yeps \to y
  & \hbox{weakly star in $\H1\Vp\cap\L\infty V\cap\L2\Hdue$}
  \qquad
  \label{convy}
  \\
  \hskip-0.8cm & \dt\yeps \to \dt y
  & \hbox{weakly in $\L2H$ \ if $\tau>0$}
  \qquad
  \label{convybis}
  \\
  \hskip-0.8cm & \yGeps \to \yG
  & \hbox{weakly star in $\H1\HG\cap\L\infty\VG\cap\L2\HdueG$}
  \label{convyG}
  \\
  \hskip-0.8cm & \weps \to w
  & \hbox{weakly in $\L2V$}
  \label{convw}
  \\
  \hskip-0.8cm & \xieps \to \xi
  & \hbox{weakly in $\L2H$}
  \label{convxi}
  \\
  \hskip-0.8cm & \xiGeps \to \xiG
  & \hbox{weakly in $\L2\HG$} 
  \label{convxiG}
\Esist
at least for a subsequence, in principle.
Clearly, $y$~satisfies the Cauchy condition \eqref{cauchy},
$\yG$~is the trace of~$y$,
and the quintuplet $(y,\yG,w,\xi,\xiG)$
satisfies \eqref{intprima} and a variational equation like \eqref{intseconda},
where the terms related to $\pi$ and $\piG$ are not yet identified.
Moreover, the relationships contained in \accorpa{regxi}{regxiG} have to be proved.
Thanks to, e.g., \cite[Sect.~8, Cor.~4]{Simon}\pier{, it is not difficult to}
infer that
\Bsist
  & \yeps \to y
  & \hbox{strongly in \,\pier{$\C0{H}$}}
  \qquad
  \label{strongy}
  \\
  & \yGeps \to \yG
  & \hbox{strongly in \,\pier{$\C0{H}$}} .
  \label{strongyG}
\Esist
\pier{Then, recalling \eqref{hppilip}}, we deduce that $\pi(\yeps)$ and $\piG(\yGeps)$
converge to $\pi(y)$ and to $\piG(\yG)$ \pier{in $\C0H$ and in $\C0\HG$}, respectively.
Moreover, by applying well-known results on maximal monotone operators
(see, e.g., \cite[Lemma~1.3, p.~42]{Barbu}),
we infer that $\xi\in\beta(y)$ \aeQ\ and $\xiG\in\betaG(\yG)$ \aeS.
This completes the proof of Theorem~\ref{Existence}.
\QED

\bigskip

The rest of the section is devoted to the proof of our regularity results.
We start with Theorem~\ref{Regularity} and thus suppose that
its assumptions are satisfied. \pier{Then we can take $\geps=g$ in both cases.}
As observed at the beginning,
just further a priori estimates on the solution to the approximating problem are necessary.
In order to confine the length of the paper,
we proceed formally, by assuming that the solution to the approximating problem 
is as smooth as needed.
We prepare a lemma.

\Blem
\label{Dtyzero}
We have
\Beq
  \normaVp{\dt\yeps(0)}
  + \taueps^{1/2} \normaH{\dt\yeps(0)}
  + \normaHG{\dt\yGeps(0)}
  \leq c \,.
  \label{dtyzero}
\Eeq
\Elem

\Bdim
The values $\dt\yeps(0)$ and $\dt\yGeps(0)$ can be obtained 
by taking $t=0$ in \accorpa{primaeps}{secondaeps}.
Hence, they satisfy
\Bsist
  \hskip-1.2cm && \< \dt\yeps(0) , v >
  + \iO \nabla\weps(0) \cdot \nabla v = 0
  \label{primaepsz}
  \\
  \noalign{\smallskip}
  \hskip-1.2cm && \iO \weps(0) v
  = \taueps \iO \dt\yeps(0) \, v
  + \iG \dt\yGeps(0) \, v
  + \iO \nabla\yz \cdot \nabla v
  + \iG \nablaG(\yz\suG) \cdot \nablaG v
  \non
  \\
  \hskip-1.2cm && \quad {}
  + \iO \bigl( \betaeps(\yz) + 
  \lam\pier{(0)}\,\pi(\yz) - \geps(0) \bigr)  v
  + \iG \bigl( \betaGeps(\yz\suG) + 
  \lamG\pier{(0)}\,\piG(\yz\suG) - \gG(0) \bigr) v \quad
  \label{secondaepsz}
\Esist
for every $v\in V$ and every $v\in\calV$, respectively.
We choose $v=\calN(\dt\yeps(0))$ and $v=-\dt\yeps(0)$ 
in \eqref{primaepsz} and in~\eqref{secondaepsz}, respectively,
sum the obtained equalities to each other 
and \pier{exploit} \eqref{defN} and~\eqref{defnormavp}, as usual.
By observing that $\dt\yGeps(0)=(\dt\yeps(0))\suG$ (since $\yeps$ is smooth),
integrating by parts both in~$\Omega$ 
(with the help of our assumption \eqref{hpyzHdue} on~$\yz$) and on~$\Gamma$, 
and rearranging a little,
we~have
\Bsist
  \hskip-1cm && \normaVp{\dt\yeps(0)}^2
  + \taueps \normaH{\dt\yeps(0)}^2
  + \normaHG{\dt\yGeps(0)}^2
  \non
  \\
  \hskip-1cm && = - \iO \bigl( -\Delta\yz + 
  \betaeps(\yz) + \lam\pier{(0)}\,\pi(\yz) - 
  \geps(0) \bigr) \dt\yeps(0)
  \non
  \\
  \hskip-1cm && \quad {}
  - \iG \bigl( -\DeltaG\yz\suG + \betaGeps(\yz\suG) + \lamG\pier{(0)}\,\piG(\yz\suG) 
  - \gG(0) \bigr) \dt\yGeps(0) .
  \non
\Esist
The last integral can be easily handled by using 
\eqref{hpyzHdue}, \eqref{hpxiGz}, \eqref{propyosida}, 
assumptions \accorpa{hppilip}{hplam} on $\piG$ and~$\lamG$, and~\eqref{hpgbis}.
\pier{Then, we deduce that}
\Beq
  - \iG \bigl( -\DeltaG\yz\suG + \betaGeps(\yz\suG) + \lamG \pier{(0)}\,
    \piG(\yz\suG) - \gG(0) \bigr) \dt\yGeps(0)
  \leq \frac 12 \, \normaH{\dt\yGeps(0)}^2
  + c \,.
  \non
\Eeq
We can deal with the integral over $\Omega$ in a similar way if $\tau>0$ 
(e.g., \eqref{hpxiz} replaces \eqref{hpxiGz} in the argument).
In~this case, we obtain
\Beq
  - \iO \bigl( -\Delta\yz + \betaeps(\yz) + \lam\pier{(0)}\,\pi(\yz) - \geps(0) \bigr) \dt\yeps(0)
  \leq \frac \taueps 2 \, \normaH{\dt\yeps(0)}^2
  + c \,.
  \non
\Eeq
On the contrary, if $\tau=0$, the treatment of the integral is more delicate 
and requires the help of \eqref{hpxizV}, \eqref{hpdxlam} and~\eqref{interpolaz}.
We~have
\begin{align}
  & - \iO \bigl( -\Delta\yz + \betaeps(\yz) + \lam\pier{(0)}\,\pi(\yz) - \geps(0) \bigr) \dt\yeps(0)
  \non
  \\
  & \leq \frac 12 \, \normaVp{\dt\yeps(0)}^2
  + c \, \normaV{-\Delta\yz + \betaeps(\yz) + \lam\pier{(0)}\,\pi(\yz) - \geps(0)}^2
  \non
  \\
  & \leq \frac 12 \, \normaVp{\dt\yeps(0)}^2
  + c \, \normaV{-\Delta\yz + \betaeps(\yz) - \geps(0)}^2
  + c \, \normaV{\lam\pier{(0)}\,\pi(\yz)}^2
  \leq \frac 12 \, \normaVp{\dt\yeps(0)}^2
  + c \,.
  \non
\end{align}
In both cases, we can combine and conclude that \eqref{dtyzero} holds true.
\Edim

\step
Sixth a priori estimate

We differentiate equations \eqref{primaeps} and \eqref{secondaeps} with respect to time
and obtain \aat\ (but we avoid writing the time~$t$ everywhere, for brevity)
\Bsist
  && \< \dt^2\yeps , v >
  + \iO \nabla\dt\weps \cdot \nabla v = 0
  \label{dtprimaeps}
  \\
  \noalign{\smallskip}
  && \iO \dt\weps v
  = \taueps \iO \dt^2\yeps \, v
  + \iG \dt^2\yGeps \, v
  + \iO \nabla\dt\yeps \cdot \nabla v
  + \iG \nablaG\dt\yGeps \cdot \nablaG v
  \qquad
  \non
  \\
  && \quad {}
  + \iO \bigl(
    \betaeps'(\yeps) \dt\yeps
    + (\dt\lam) \, \pi(\yeps) + \lam \, \pi'(\yeps) \dt\yeps - \dt\geps
  \bigr) \, v
  \non
  \\
  && \quad {}
  + \iG \bigl(
    \betaGeps'(\yGeps) \dt\yGeps
    + (\dt\lamG) \piG(\yGeps) + \lamG \piG'(\yGeps) \dt\yGeps - \dt\gG
  \bigr) \, v
  \label{dtsecondaeps}
\Esist
where \eqref{primaeps} and \eqref{secondaeps} are required to hold
for every $v\in V$ and every $v\in\calV$, respectively.
Now, we note that \eqref{conservedeps} implies $\dt\yeps(t)\in D(\calN)$ \aat\
(cf.~\eqref{predefN}).
So, we test the above equations by $\calN(\dt\yeps)$ and~$-\dt\yeps$,
integrate over~$(0,t)$, sum up, and account for \eqref{defN} and~\eqref{defnormavp}, as usual.
We obtain
\Bsist
  && \frac 12 \, \normaVp{\dt\yeps(t)}^2
  + \frac \taueps 2 \iO |\dt\yeps(t)|^2
  + \frac 12 \iG |\dt\yGeps(t)|^2
  \non
  \\
  && \quad
  + \intQt |\nabla\dt\yeps|^2
  + \intSt |\nablaG\dt\yGeps|^2
  + \intQt \betaeps'(\yeps) |\dt\yeps|^2
  + \intSt \betaGeps'(\yGeps) |\dt\yGeps|^2
  \non
  \\
  && = \frac 12 \,\normaVp{\dt\yeps(0)}^2
  + \frac \taueps 2 \iO |\dt\yeps(0)|^2
  + \frac 12 \iG |\dt\yGeps(0)|^2
  \non
  \\
  && \quad
  + \intQt \bigl( \dt\geps - (\dt\lam) \pi(\yeps) - \lam \, \pi'(\yeps) \dt\yeps \bigr) \dt\yeps
  \non
  \\
  && \quad
  + \intSt \bigl( \dt\gG - (\dt\lamG)\piG(\yGeps) - \lamG \piG'(\yGeps) \dt\yGeps \bigr) \dt\yGeps \,.
  \label{persestastima}
\Esist
All the integrals on the \lhs\ are nonnegative
and the first three terms on the \rhs\ are bounded thanks to Lemma~\ref{Dtyzero}.
The next integral is estimated as follows
\Bsist
  && \intQt \bigl( \dt\geps - (\dt\lam) \pi(\yeps) - \lam \, \pi'(\yeps) \dt\yeps \bigr) \dt\yeps
  \non
  \\
  && \leq c \intQt |\dt\yeps|^2
  + \norma{\dt\geps}_{\L2H}^2
  + \norma{(\dt\lam)\pi(\yeps)}_{\L2H}^2
  \non
  \\
  && \leq c \intQt |\dt\yeps|^2
  + c
  + c \norma{\dt\lam}_{\L\infty H}^2 \norma\yeps_{\L2\Linfty}^2
  \non
  \\
  && \leq \frac 12 \intQt |\nabla\dt\yeps|^2
  + c \iot \normaVp{\dt\yeps(s)}^2 \, ds
  + c
  \non
\Esist
thanks to \eqref{hpdtlam}, the boundedness of~$\pi'$ (cf.~\eqref{hppilip}), \eqref{hpgbis},
the interpolation inequality~\eqref{interpolaz} and~\eqref{dastimeHdue}.
As the last integral can be treated in a similar and even simpler way,
\pier{from \eqref{persestastima}} we conclude~that
\Bsist
  && \norma{\dt\yeps}_{\L\infty\Vp}
  + \taueps^{1/2} \norma{\dt\yeps}_{\L\infty H}
  + \norma{\dt\yGeps}_{\L\infty\HG}
  \non
  \\
  && \quad {}
  + \norma{\dt\yeps}_{\L2V}
  + \norma{\dt\yGeps}_{\L2\VG} 
  \leq c \,.
  \label{sestastima}
\Esist

\step
Conclusion of the proof of Theorem~\ref{Regularity}

Due to the existence proof, 
we already know that $(\yeps,\yGeps,\weps,\xieps,\xiGeps)$
converges to the solution to problem \Pbl.
Thus, the last estimate of ours just improves 
the regularity of the limit
(as~well as the topology of the convergence),
and \accorpa{regybis}{regyGbis} are partially proved.
In order to achieve the $H^2$ regularity requirements of the statement,
one can argue as in \pier{the proof of Theorem~\ref{Existence}, going from the Third a priori estimate to the Conclusion of the proof.}
In the argument used there, time was just a parameter, indeed,
since everything was based on the theory of elliptic regularity and traces in~$\Omega$.
In the present case one immediately checks that \pier{$L^\infty$ instead of $L^2$ bounds hold} with respect to time.
Hence, we \pier{get, in this order, 
\begin{align}
  &\norma\weps_{\L\infty V} \leq c \label{pier1}\\
  &\norma\xieps_{\L\infty H} \leq c \label{pier2}\\
  &\norma\xiGeps_{\L\infty H} \leq c \label{pier3}\\
  &\norma\yeps_{\L\infty\Hdue}
  + \norma\yGeps_{\L\infty\HdueG} \leq c \label{stimayHdue}
\end{align}
}%
i.e., the last part of \accorpa{regybis}{regyGbis} \pier{as well as \eqref{pier0}}. This completes the proof.
\QED

\step
Proof of Theorem \ref{Bddness}

We recall that $\tau>0$.
\pier{Hence}, $\dt\yeps$ is bounded in~$\L\infty H$.
On the other hand, equation \eqref{primaeps} implies that
$\dt\yeps-\Delta\weps=0$ in~$Q$ and $\dn\weps=0$ on the boundary, 
whence 
\Beq
  \norma{\Delta\weps}_{\L\infty H} \leq c 
  \aand
  \norma\weps_{\L\infty\Hdue} \leq c \,.
  \label{stimawHdue}
\Eeq
This implies the first \pier{conclusion in} \eqref{bddness}.
In order to prove boundedness for~$\xi$,
it suffices to find an a~priori estimate for the $L^p$ norm of~$\xieps$ 
that is uniform with respect to both $p$ and~$\eps$.
Thus, in the sequel, the dependence on $p$ of the constants 
is explicitly written and carefully controlled.
In order to perform our estimate, we write \eqref{secondaeps}~as
\begin{align}
  & \tau \iO \dt\yeps \, v
  + \iG \dt\yGeps \, v
  + \iO \nabla\yeps \cdot \nabla v
  + \iG \nablaG\yGeps \cdot \nablaG v
  + \iO \betaeps(\yeps) v 
  + \iG \betaGeps(\yGeps) v 
  \non
  \\
  & = \iO \feps v + \iG \fGeps v
  \quad \hbox{\aet\ and for every $v\in\calV$}
  \label{perbdd}
\end{align}
where $\feps:=w+g-\lam\,\pi(\yeps)$ and $\fGeps:=\gG-\lamG\piG(\yGeps)$,
and observe that 
\Beq
  \norma\feps_{\LQ\infty} \leq c
  \aand
  \norma\fGeps_{\LS\infty} \leq c
  \label{rhsbdd}
\Eeq
due to \eqref{hpbdd} and \accorpa{stimayHdue}{stimawHdue}.
Then, we test \eqref{perbdd} by $|\betaeps(\yeps)|^{p-1}\sign\yeps$
with an arbitrary $p>2$, where the $\sign$ function is extended by $\sign0=0$,
and integrate over~$(0,t)$.
We~have
\begin{align}
  & \tau \iO \Bp(\yeps(t))
  + \iG \Bp(\yGeps(t))
  \non
  \\
  & \quad {}
  + (p-1) \iO \betaeps'(\yeps) |\betaeps(\yeps)|^{p-2} |\nabla\yeps|^2
  + (p-1) \iG \betaeps'(\yGeps) |\betaeps(\yGeps)|^{p-2} |\nablaG\yGeps|^2
  \non
  \\
  & \quad {}
  + \intQt |\betaeps(\yeps)|^p
  + \intSt \betaGeps(\yGeps) \, |\betaeps(\yeps)|^{p-1}\sign\yGeps
  \non
  \\
  & = \tau \iO \Bp(\yz)
  + \iG \Bp(\yz\suG)
  \non
  \\
  & \quad {}
  + \intQt \feps \,|\betaeps(\yeps)|^{p-1}\sign\yeps
  + \intSt \fGeps |\betaeps(\yGeps)|^{p-1}\sign\yGeps 
  \label{perxibdd}
\end{align}
where we have set
\Beq
  \Bp(r) := \int_0^r |\betaeps(s)|^{p-1} \sign s \, ds
  \quad \hbox{for $r\in\erre$} .
  \label{defBp}
\Eeq
We recall that $\betaeps$ and $\betaGeps$ are monotone functions that vanish at the origin.
It follows that they and the identity map have the same sign,
whence all the terms on the \lhs\ of \eqref{perxibdd} are nonnegative.
Moreover, the last of them can be estimated from below 
on account of the compatibility condition~\eqref{hpcompatib} 
and of the Young inequality~\eqref{young}
(with $p'$ in place of~$p$ and $\delta>0$ to be chosen) as~follows
\Bsist
  && \intSt \betaGeps(\yGeps) \, |\betaeps(\yeps)|^{p-1}\sign\yGeps
  \geq \intSt \bigl( \eta |\betaeps(\yGeps)|^p - C |\betaeps(\yGeps)|^{p-1} \bigr)
  \non
  \\
  && \geq \eta \intSt |\betaeps(\yGeps)|^p
  - \intSt \Bigl( \frac \delta {p'} |\betaeps(\yGeps)|^{(p-1)p'} + \frac {\delta^{-p/p'}} p \, C^p \Bigr) \non
  \\
  && = \eta \intSt |\betaeps(\yGeps)|^p
  - \frac \delta {p'} \intSt |\betaeps(\yGeps)|^p
  - c \, \frac {\delta^{-p/p'}} p \, C^p .
  \non
\Esist
By choosing $\delta=\eta p'/2$, we conclude that
\Beq
  \intSt \betaGeps(\yGeps) \, |\betaeps(\yeps)|^{p-1}\sign\yGeps
  \geq \frac \eta 2 \intSt |\betaeps(\yGeps)|^p
  - c^p \,.
  \label{frombelow}
\Eeq
Let us come to the \rhs\ and denote by $M$ the $L^\infty$ norm of~$\betaz(\yz)$ (cf.~\eqref{hpbdd}).
By~\eqref{propyosida}, we deduce that
\Beq
  |\betaeps(\yz)| \leq M
  \quad \aeO
  \aand
  |\betaeps(\yz\suG)| \leq M
  \quad \aeG .
  \non
\Eeq
Thus, on account of \eqref{defBp} and of the $L^\infty$ bound that follows from~\eqref{stimayHdue},
we can estimate the sum of the first two integrals this~way
\Beq
  \tau \iO \Bp(\yz)
  + \iG \Bp(\yz\suG)
  \leq c M^{p-1} 
  \leq c^p .
  \non
\Eeq
Now, we consider the volume integral.
We recall \eqref{rhsbdd}, apply the Young inequality and have
\Beq
  \iO \feps \,|\betaeps(\yeps)|^{p-1}\sign\yeps
  \leq \intQt \Bigl( \frac 1p \, c^p + \frac 1{p'} |\betaeps(\yeps)|^{(p-1)p'} \Bigr) 
  \leq c^p + \frac 1{p'} \intQt |\betaeps(\yeps)|^p .
  \non
\Eeq
Arguing as for \eqref{frombelow}, we obtain
\Beq
  \intSt \fGeps |\betaeps(\yGeps)|^{p-1}\sign\yGeps 
  \leq \frac \eta 4 \intSt |\betaeps(\yGeps)|^p 
  + c^p .
  \non
\Eeq
By collecting \eqref{perxibdd} (where we neglect a number of nonnegative terms on the \lhs), 
\eqref{frombelow} and the last two estimates, and rearranging,
we infer~that
\Beq
  \frac 1p \intQt |\betaeps(\yeps)|^p 
  + \frac \eta 4 \intSt |\betaeps(\yGeps)|^p
  \leq c^p
  \non
\Eeq
and easily conclude that
\Beq
  \norma{\betaeps(\yeps)}_{\LQ p}
  + \norma{\betaeps(\yGeps)}_{\LS p}
  \leq c \,.
  \non
\Eeq
This completes the proof.
\QED

\Brem
The above proof also provides an $L^\infty$ bound for~$\betaeps(\yGeps)$.
This implies that
\Beq
  \betaeps(\yGeps) \to \zeta
  \quad \hbox{weakly star in $\LS\infty$}
  \non
\Eeq
for some $\zeta\in\LS\infty$ and at least for a subsequence.
This and the strong convergence of $\yGeps$, e.g., in $\LS2$ yield
$\zeta\in\beta(\yG)$ by maximal monotonicity.
Hence, we have also proved that some selection of $\beta(\yG)$ is bounded on~$\Sigma$.
On the contrary, nothing can be inferred as far as $\betaG(\yG)$ is concerned,
unless the assumptions of Corollary~\ref{StrongBddness} are supposed to hold.
\Erem

%%%%%%%%%%%%%%%%%%%%%%%%%%%%%%%%%%%%%%%%%%%%%%%%%%%%%%%%%%%%%%%%%%%%%%%%

%%%%%%%%%%%%%%%%%%%%%%%%%%%%%%%%%
%% bibliography
%%%%%%%%%%%%%%%%%%%%%%%%%%%%%%%%%

\vspace{3truemm}

\Begin{thebibliography}{10}

\bibitem{Barbu}
V. Barbu,
``Nonlinear semigroups and differential equations in Banach spaces'',
Noord-hoff, 
Leyden, 
1976.

\bibitem{Brezis}
H. Brezis,
``Op\'erateurs maximaux monotones et semi-groupes de contractions
dans les espaces de Hilbert'',
North-Holland Math. Stud.
{\bf 5},
North-Holland,
Amsterdam,
1973.

\bibitem{BreGil}
F. Brezzi and G. Gilardi,
Part~1:
Chapt.~1, Functional spaces, 
Chapt.~2, Partial differential equations,
in ``Finite element handbook'',
H. Kardestuncer and D.H. Norrie eds.,
McGraw-Hill Book Company,
NewYork,
1987.

\bibitem{CahH} 
J.W. Cahn and J.E. Hilliard, 
Free energy of a nonuniform system I. Interfacial free energy, 
{\it J. Chem. Phys.\/}
{\bf 2} (1958) 258-267.

\bibitem{CaCo}
L. Calatroni and P. Colli,
Global solution to the Allen-Cahn equation with singular potentials and dynamic boundary conditions,
{\it \pier{Nonlinear Anal.}\/} 
{\bf 79} (2013) 12-27.

\bibitem{CFP} 
R. Chill, E. Fa\v sangov\'a, and J. Pr\"uss,
Convergence to steady states of solutions of the Cahn-Hilliard equation with dynamic boundary conditions,
{\it Math. Nachr.\/} 
{\bf 279} (2006) 1448-1462.

\pier{
\bibitem{CFS}
P. Colli, M.H. Farshbaf-Shaker, and J. Sprekels,
A deep quench approach to the optimal control of an Allen-Cahn equation
with dynamic boundary conditions and double obstacles,
{\it Appl. Math. Optim.\/}, in print (2014) {\tt doi:10.1007/s00245-014-9250-8}}

\pier{
\bibitem{CGPS3} 
P. Colli, G. Gilardi, P. Podio-Guidugli, and J. Sprekels,
Well-posedness and long-time behaviour for 
a nonstandard viscous Cahn-Hilliard system, 
{\it SIAM J. Appl. Math.} {\bf 71} (2011) 1849-1870.
}

\pier{
\bibitem{CS}
P. Colli and J. Sprekels,
Optimal control of an Allen-Cahn equation with singular potentials and dynamic boundary condition,
preprint arXiv:1212.2359~[math.AP] (2012), pp.~1-24.
}%

\bibitem{EllSt} 
C.M. Elliott and A.M. Stuart, 
Viscous Cahn-Hilliard equation. II. Analysis, 
{\it J. \pier{Differential Equations}\/} 
{\bf 128} (1996) 387-414.

\bibitem{EZ} 
C.M. Elliott and S. Zheng, 
On the Cahn-Hilliard equation, 
{\it Arch. Rational Mech. Anal.\/} 
{\bf 96} (1986) 339-357.

\pier{
\bibitem{FG} 
E. Fried and M.E. Gurtin, 
Continuum theory of thermally induced phase transitions based on an order 
parameter, {\it Phys. D} {\bf 68} (1993) 326-343.
}

\bibitem{GiMiSchi} 
G. Gilardi, A. Miranville, and G. Schimperna,
On the Cahn-Hilliard equation with irregular potentials and dynamic boundary conditions,
{\it Commun. Pure Appl. Anal.\/} 
{\bf 8} (2009) 881-912.

\bibitem{GiMiSchi2} 
G. Gilardi, A. Miranville, and G. Schimperna,
Long-time behavior of the Cahn-Hilliard equation with irregular potentials and dynamic boundary conditions,
{\it \pier{Chin. Ann. Math. Ser. B}\/} 
{\bf 31} (2010) 679-712. 

\bibitem{Gu} M. Gurtin, 
Generalized Ginzburg-Landau and Cahn-Hilliard equations based on a microforce balance,
{\it Phys.~D\/} {\bf 92} (1996) 178-192.

\bibitem{Is} H. Israel,
Long time behavior of an {A}llen-{C}ahn type equation with a
singular potential and dynamic boundary conditions,
{\it J. Appl. Anal. Comput.\/} {\bf 2} (2012), 29-56.

\bibitem{LioMag}
J.-L. Lions and E. Magenes,
``Non-homogeneous boundary value problems and applications'',
Vol.~I,
Springer, Berlin, 1972.

\bibitem{MirZelik} 
A. Miranville and S. Zelik,
Robust exponential attractors for Cahn-Hilliard type equations with singular potentials,
{\it Math. \pier{Methods} Appl. Sci.\/} 
{\bf 27} (2004) 545--582.

\pier{
\bibitem{Podio}
P. Podio-Guidugli, 
Models of phase segregation and diffusion of atomic species on a lattice,
{\it Ric. Mat.} {\bf 55} (2006) 105-118.
}

\bibitem{PRZ} 
J. Pr\"uss, R. Racke, and S. Zheng, 
Maximal regularity and asymptotic behavior of solutions for the Cahn-Hilliard equation with dynamic boundary conditions,  
{\it Ann. Mat. Pura Appl.~(4)\/}
{\bf 185} (2006) 627-648.

\bibitem{RZ} 
R. Racke and S. Zheng, 
The Cahn-Hilliard equation with dynamic boundary conditions, 
{\it Adv. \pier{Differential Equations}\/} 
{\bf 8} (2003) 83-110.

\bibitem{Simon}
J. Simon,
{Compact sets in the space $L^p(0,T; B)$},
{\it Ann. Mat. Pura Appl.~\pier{(4)}\/} 
{\bf 146} (1987) 65-96.

\bibitem{WZ} H. Wu and S. Zheng.,
Convergence to equilibrium for the Cahn-Hilliard equation with dynamic boundary conditions, 
{\it J. \pier{Differential Equations}\/}
{\bf 204} (2004) 511-531.

\End{thebibliography}

\End{document}

\bye